\documentclass[10pt]{article}
\usepackage{amssymb,amsmath}
\usepackage{amsthm}
\usepackage{graphicx}
\usepackage{enumerate}
\usepackage{algorithmic}
\usepackage{booktabs}

\def\vec#1{{\boldsymbol #1}}

\def\ds{\displaystyle}

\newcommand{\Dx}{\partial_{\boldsymbol{x}}}

\def\pfrac#1/#2.{\frac{\partial #1}{\partial #2}}

\def\ne{n_{\rm e}}
\def\np{n_{\rm p}}
\def\nn{n_{\rm n}}
\def\ve{\vec{v}_{\! \rm e}}
\def\vp{\vec{v}_{\! \rm p}}
\def\vn{\vec{v}_{\! \rm n}}
\def\De{D_{\rm e}}
\def\Dp{D_{\rm p}}
\def\Dn{D_{\rm n}}
\def\Sph{S_{\rm ph}}

\def\pOtwo{p_{\rm O_2}}
\def\light{c}

\newcommand{\R}{{\mathbb R}}

\newcommand{\Z}{{\mathbb Z}}

\newcommand{\M}{{\mathcal M}_n(\R)}

\def\ds{\displaystyle}

\newcommand{\Or}{{\mathcal O}}
\newcommand{\Time}{{\mathcal T}}

\newcommand{\Id}{{\rm Id}}
\newcommand{\der}{{\mathrm{d}}}

\def\infi{\infty}
\newcommand{\x}{{\boldsymbol{x}}}

\newcommand{\MR}{{\mathcal M}}
\newcommand{\thr}{{\mathcal T}}
\newcommand{\adap}{{\mathcal A}}
\newcommand{\leaf}{{\mathrm{L}}}
\newcommand{\tree}{{\mathrm{T}}}
\newcommand{\phan}{{\mathrm{P}}}

\def\disc#1{{\mathbf #1}}

\def\hypre{{\sl hypre} {}}

\title{A numerical strategy to discretize and solve the Poisson equation on 
dynamically adapted multiresolution grids
for time-dependent streamer discharge simulations}
\author{
Max~Duarte\footnotemark[2]~\footnotemark[3]
\and Zden\v{e}k~Bonaventura\footnotemark[4]
\and Marc~Massot\footnotemark[5]~\,\footnotemark[6]~\footnotemark[7]
\and Anne~Bourdon\footnotemark[5]~\,\footnotemark[6]
}

\begin{document}

\maketitle

\renewcommand{\thefootnote}{\fnsymbol{footnote}}

\footnotetext[2]{
Univ. Nice Sophia Antipolis, CNRS,  LJAD, UMR 7351, 06100 Nice, France.
}

\footnotetext[3]{
CCSE,
Lawrence Berkeley National Laboratory,
1 Cyclotron Rd. MS 50A-1148,
94720 Berkeley, CA, USA
({\tt MDGonzalez@lbl.gov}).
}

\footnotetext[4]{
Department of Physical Electronics, Faculty of Science, Masaryk
University, Kotl\'a\v{r}sk\'a 2, 61137 Brno, Czech Republic
({\tt zbona@physics.muni.cz}).
}

\footnotetext[5]{
CNRS, UPR 288, Laboratoire EM2C, 
Grande voie des vignes, 92295
Ch\^{a}tenay-Malabry, France.}

\footnotetext[6]{
Ecole Centrale Paris, 92295 Ch\^{a}tenay-Malabry,
France ({\tt \{marc.massot,anne.bourdon\}@ecp.fr}).
}

\footnotetext[7]{
F\'ed\'eration de Math\'ematiques de l'Ecole Centrale Paris, FR CNRS 3487.}

\renewcommand{\thefootnote}{\arabic{footnote}}

\begin{abstract}
We develop a numerical strategy to solve 
multi-dimensional Poisson equations
on dynamically adapted grids 
for evolutionary 
problems
disclosing propagating fronts.
The method is an extension of the multiresolution finite volume scheme
used to solve hyperbolic and parabolic time-dependent PDEs.
Such an approach
guarantees a numerical 
solution
of the Poisson 
equation within a user-defined
accuracy tolerance.
Most adaptive meshing approaches 
in the literature 
solve 
elliptic PDEs
level-wise and hence at uniform resolution throughout the set of adapted grids.
Here
we introduce a numerical procedure
to represent 
the elliptic operators on the adapted grid,
strongly coupling inter-grid relations
that guarantee the conservation and accuracy properties
of multiresolution 
finite volume schemes.
The discrete Poisson equation
is solved at once over the entire computational domain
as a completely separate process.
The accuracy and numerical performance of the method 
is assessed
in the context of
streamer discharge
simulations.
\end{abstract}

\subsubsection*{Keywords} 
Poisson equation,  
multiresolution finite volume scheme, streamer discharges.

\subsubsection*{AMS subject classifications}
35J05, 65M50, 65G20, 65M08, 65Y15, 65Z05, 76X05

\pagestyle{myheadings}
\thispagestyle{plain}
\markboth{DUARTE, BONAVENTURA, MASSOT, BOURDON}
{POISSON EQUATION ON ADAPTIVE MULTIRESOLUTION GRID}

\section{Introduction}\label{sec:intro}
In numerous scientific applications
one has to deal with the numerical
solution of elliptic PDEs, like Poisson equations,
coupled with evolutionary PDEs 
to address the numerical simulation of time-dependent physical processes.
One major example is given, for instance, 
by the so-called projection methods \cite{Chorin68,Temam2}, 
widely investigated, extended, 
and implemented in the literature
to solve the incompressible Navier-Stokes equations
(see, {\it e.g.}, \cite{Guermond2006} and references therein).
Solving Poisson equations is also very common in plasma physics simulations.
As an example, in the framework of a drift-diffusion model consisting of a set 
of continuity equations for charged species coupled with a Poisson equation for the electric potential,  
non-linear ionization waves also called streamers can be simulated \cite{Babaeva:1996,Kulikovsky:1997c}. 
In either situation
Poisson-type equations
must be solved (often several times) at every time-step
throughout the numerical simulation,
a task that
depending on the size and complexity
of the problem can easily 
become cumbersome in both CPU time and memory.
In particular
phenomena characterized by propagating fronts,
as considered in this work, 
commonly require a sufficiently fine spatial
representation and potentially large
systems of equations need then to be solved.

In this regard
grid adaptation for  time-dependent problems
disclosing localized fronts 
is specifically designed to yield high data compression and hence
important savings in computational costs
(see, {\it e.g.}, \cite{berger1984,berger1989}).
Among the many adaptive meshing approaches developed in the literature,
we consider in this work adaptive multiresolution
schemes based on \cite{Harten94,Harten95},
namely the multiresolution finite volume
scheme introduced in \cite{Cohen03}
for conservation laws.
Besides the inherent advantages
of grid adaptation, multiresolution
techniques rely on biorthogonal wavelet decomposition \cite{Cohen1992}
and thus offer a rigorous 
mathematical framework for adaptive meshing
schemes \cite{cohen2000a,muller2003}.
In this way
not only approximation errors 
coming from grid adaptation and thus data compression
can be tracked, but general and robust 
techniques can be built since 
the wavelet decomposition is
independent of any physical particularity 
of the problem and accounts
only for the spatial regularity of the 
discrete variables at a given simulation time.

Adaptive multiresolution schemes have been successfully
implemented for the simulation of compressible fluids modeled
by Euler or Navier-Stokes equations
(see, {\it e.g.}, \cite{Muller2009,Brix2011,Domingues2011} 
and references therein),
as well as for the numerical solution of time-dependent
parabolic \cite{Roussel03,Burger08} and stiff parabolic PDEs 
\cite{Duarte11_SISC,Dumont2013,DuarteCFlame}.
Nevertheless, to the best of our knowledge this is the fist attempt
to develop a Poisson solver in the context of the adaptive multiresolution 
finite volume method introduced in \cite{Cohen03}
for evolutionary problems.
Previously, such a solver was introduced in
\cite{Vasilyev2005} in the the context of wavelet collocation methods 
for evolutionary PDEs developed in \cite{Vasilyev2000,Vasilyev2003}.
Analogous to multiresolution schemes, wavelet collocation methods
assure adaptive meshing capabilities within a user-defined accuracy
exploiting the mathematical properties of wavelet decomposition
(for a recent review on wavelet methods see \cite{Schneider2010} 
and references therein).
Notice that an important amount of research 
has been conducted in the past decades
to solve elliptic PDEs using wavelet methods and multiresolution
representations 
(see, for instance, \cite{Jaffard92,Dahlke1997,CohenMasson1999,CohenMasson2000,Barinka2001}).
In this context the numerical solution of an elliptic PDE
is in general performed
using compressed representations of the problem in an
appropriate wavelet space, within a solid mathematical framework 
(see \cite{Cohen2001} and references therein).
Here we do not consider wavelet methods to solve elliptic PDEs,
but rather aim at developing a numerical strategy
to discretize and solve Poisson equations on dynamically adapted finite volume grids
generated by means of a multiresolution analysis.

Dynamic meshing techniques for finite volume 
discretizations are usually implemented 
by defining a set of embedded grids with different spatial resolution.
Particular attention must be addressed to 
the inter-grid interfaces in order to 
consistently define the discrete operations there.
Otherwise, potential mismatches 
may lead to substantial differences in the numerical
approximations as well as loss of conservation
(see \cite{Almgren1998} for a detailed discussion).
The most common way of solving an elliptic PDE 
on this type of adapted grid
consists in solving the 
discrete system level-wise, that 
is, considering one grid-level at a time
followed by inter-level operations to 
synchronize shared interfaces at different
grid-levels 
as well as overlapped regions.
Computations are thus successively performed
over partial regions 
at
a uniform mesh resolution until the
problem is entirely solved on the adapted grid.
Some examples can be found, for instance, 
in \cite{Almgren1998,Teyssier2002,Montijn:2006,Martin2008,Safta2010}.
For intensive computations
iterative linear solvers based on
geometric multigrid schemes are often implemented, 
taking advantage of
the multi-mesh representation of the problem 
\cite{Almgren1998,Teyssier2002,Martin2008}.
In particular
the Poisson solver in \cite{Vasilyev2005}
also implements a level-wise approach where a finite difference
discretization is considered.

The main objective of this paper is to develop a 
Poisson solver on dynamically adapted grids
generated with a multiresolution finite volume scheme.
In particular
we investigate
the influence of data compression on the
accuracy of approximations
obtained with Poisson equations discretized
on an adapted multiresolution mesh.
One novelty of this paper in terms of elliptic solvers
on adapted grids is that
instead of solving the discrete equations
level-wise throughout the set of embedded
grids, we have conceived a numerical procedure
to represent the elliptic operators
discretized directly on the adapted grid,
that is, on a mesh
consisting of cells with different spatial resolution.
The algorithm relies on a local
reconstruction of 
uniform-grid zones at inter-level interfaces
by means of multiresolution operations
between consecutive grid-levels
that guarantee the conservation and accuracy properties
of multiresolution schemes.
This 
approach 
results in a separate
algebraic system completely
independent of any consideration
related to the adaptive meshing scheme
or its corresponding data structure, as well as 
of the numerical integration of 
the time-dependent PDEs
associated with the model.
The resulting
discrete system
can thus be solved 
at once over the whole computational domain
with no need of grid overlapping
by considering an appropriate linear solver.

The performance of the strategy
is assessed
in the context of streamer discharge simulations at atmospheric pressure.
The detailed physics of these discharges reveals an important 
time-space multi-scale character \cite{Ebert_nonlinearity:2011}.
Grid adaptation is therefore highly desirable and 
was already considered, 
for instance, in \cite{Montijn:2006,Pancheshnyi:2008,Unfer:2010}.
In \cite{Duarte11_JCP}
we introduced
a time-space adaptive numerical scheme 
with error control
to simulate
propagating streamers on multiresolution
grids.
Nevertheless, 
a simplified geometry was considered there in order to
avoid the numerical solution of a
multi-dimensional
Poisson equation.
The present work describes the required fundamentals and 
further developments needed to 
solve Poisson equations on a finite volume adapted grid
according to the approach established in \cite{Duarte11_JCP}.
The latter aims at assuring 
a tracking capability of the numerical errors 
and a 
full resolution of the equations on the adapted grid.

The paper is organized as follows.
We give in Section \ref{sec:PoissonMR}
a short introduction on multiresolution finite volume schemes
and describe the data compression errors
associated with 
Poisson equations discretized
on multiresolution grids.
In Section \ref{sec:imple}
we recall 
some key aspects of the
multiresolution technique
considered here.
We then describe the numerical procedure
conceived to represent elliptic operators on the finite volume adapted mesh.
Numerical results coming from 
streamer discharge simulations are
investigated in Section \ref{sec:streamer}.

\section{Data compression errors for Poisson equations on multiresolution grids}\label{sec:PoissonMR}
We investigate 
the impact of data compression on the numerical accuracy of the approximations
obtained with a Poisson equation discretized on a multiresolution adapted grid.
However, we first need to 
briefly recall the general framework of multiresolution finite volume schemes.
More details on wavelet decomposition and multiresolution
techniques for grid adaptation
can be found in \cite{cohen2000a,muller2003}.

\subsection{Multiresolution finite volume scheme}
According to the multiresolution finite volume scheme
\cite{Cohen03},
let us build a set of nested dyadic grids 
over $\Omega \subset \mathbb{R}^{d}$ as follows.
We consider regular disjoint partitions (cells)
$(\Omega_{\gamma})_{\gamma \in S_j}$
of $\Omega$ such that
$\bigcup_{\gamma \in S_j} \Omega_{\gamma} = \Omega$
for $j=0,1,\ldots,J$.
Since each $\Omega _{\gamma}$, $\gamma \in S_j$, 
is the union of a finite number of cells $\Omega_{\mu}$
($2^d$ cells in the dyadic case), $\mu \in S_{j+1}$,
the sets $S_j$ and $S_{j+1}$ represent consecutive embedded grids over $\Omega$,
where $j$ corresponds to the grid-level from
the coarsest $(j=0)$ to the finest $(j=J)$ grid.
Defining $\Omega_{\gamma}:=\Omega_{j,k}$, 
we denote $|\gamma|:=j$ if $\gamma \in S_j$, while subscript $k\in \Delta_j \subset \Z^d$
corresponds to the position of the cell within $S_j$. 
For instance, 
in Cartesian coordinates
we consider the 
univariate dyadic intervals in $\R$:
\begin{equation}\label{eq3:dyadic_1D}
\Omega_{\gamma}=
\Omega_{j,k} := ]2^{-j}k,2^{-j}(k+1)[,\ \gamma \in S_j:= \{(j,k) \ \mathrm{s.t.} \ j\in (0,1,\ldots,J), \, k\in \Z \},
\end{equation}
and the same follows for higher dimensions.

We denote
$\disc f_j:=(f_{\gamma})_{\gamma \in S_j}$ as the 
spatial representation of $f$
on the grid $S_j$, where $f_{\gamma}$
represents the cell-average of 
$f :\, \R \times \R^d \to \R$
in $\Omega _{\gamma}$:
\begin{equation}\label{eq3:average_finite_vol}
  f_{\gamma} := |\Omega_{\gamma}|^{-1} \int_{\Omega_{\gamma}}  f(t,\x)\, \der \x, \quad \x\in \R^d.
 \end{equation}
Data at different levels of discretization are related
by two inter-level transformations which are defined as follows.
First, the
{\it projection} operator $P^j_{j-1}$
maps $\disc f_j$ to $\disc f_{j-1}$.
It is obtained
through exact a\-ve\-ra\-ges computed at the finer level by
\begin{equation}\label{eq3:projection}
  f_{\gamma} = |\Omega_{\gamma}|^{-1} 
\sum_{|\mu|=|\gamma| +1,\Omega_{\mu} 
\subset \Omega_{\gamma}} |\Omega_{\mu}|f_{\mu}.
 \end{equation}
As far as grids are nested, this
projection operator is {\it exact} and {\it unique}
\cite{cohen2000a}.
Second,
the {\it prediction} operator $P^{j-1}_j$
maps $\disc f_{j-1}$ to 
an approximation $\widehat{\disc f}_j$ of $\disc f_{j}$.
Here
a polynomial interpolation
of order $\beta$ is used to define the prediction operator:
\begin{equation}\label{eq3:linear_prediction}
  \widehat{f}_{\mu} =
 \sum_{\gamma \in R_I(\mu)}
\beta_{\mu,\gamma} f_\gamma, \quad
|\mu|=|\gamma| +1,
  \end{equation}
for a set of coefficients $(\beta_{\mu,\gamma})_{\gamma \in R_I(\mu)}$
and an {\it interpolation stencil} $R_I(\mu)$
surrounding $\Omega_{\mu}$ at the coarser level $|\gamma|=|\mu|-1$.
In particular
the prediction must be {\it consistent} with the projection \cite{Cohen03} 
in the sense that
 \begin{equation}\label{eq3:prediction_consist}
   f_{\gamma} = 
   |\Omega_{\gamma}|^{-1}\sum_{|\mu|=|\gamma| +1,\Omega_{\mu} \subset \Omega_{\gamma}} |\Omega_{\mu}|\widehat{f}_{\mu};
 \end{equation}
    \textit{i.e.}, one can retrieve the coarse cell-averages from the
predicted values: 
$P_{j-1}^j \circ P_j^{j-1} = \Id$.

With these operators we define for each cell $\Omega _{\mu}$
the prediction error or {\it detail} 
as the difference between the exact and predicted values,
\begin{equation}\label{eq3:MR_detail}
d_{\mu} := u_{\mu} - \widehat{u}_{\mu},
\end{equation}
or in terms of inter-level operations:
$d_{\mu} = u_{\mu} - P_{|\mu|}^{|\mu|-1} \circ P_{|\mu|-1}^{|\mu|} u_{\mu}$.
The consistency property (\ref{eq3:prediction_consist}) and the 
definitions of the 
projection operator (\ref{eq3:projection}) and 
of the detail (\ref{eq3:MR_detail}) imply that
\begin{equation}\label{eq3:MR_consis_2}
\sum_{|\mu|=|\gamma| +1,\Omega_{\mu} \subset \Omega_{\gamma}} d_{\mu}=0.
 \end{equation}
We can then construct as shown in \cite{Cohen03}
a {\it detail vector} defined as $\disc d_j:=(d_{\mu})_{\mu \in \nabla_j}$,
where the set $\nabla_j \subset S_j$ is obtained by removing for
each $\gamma \in S_{j-1}$ one $\mu \in S_j$ ($\Omega_{\mu} \subset \Omega_{\gamma}$)
in order to avoid redundancy (considering (\ref{eq3:MR_consis_2}))
and to get a one-to-one correspondence:
\begin{equation*}\label{eq3:one_one_cor}
\mathbf{f}_{j+1}\longleftrightarrow (\mathbf{f}_j,\mathbf{d}_{j+1}),
\end{equation*}
that is, $\mathbf{f}_{j+1}$ can be exactly computed using the cell-averages
$\mathbf{f}_{j}$ at a coarser level and the set of details $\mathbf{d}_{j+1}$
defined with operators $P^j_{j-1}$ and $P_j^{j-1}$.
By iterating this decomposition, we finally obtain a 
multi-scale representation of $\mathbf{f}_J$
in terms of $\mathbf{m}_J := (\mathbf{f}_0,\mathbf{d}_1,\mathbf{d}_2,\cdots,\mathbf{d}_J)$:
\begin{equation}\label{eq3:M_U_M}
\MR:\mathbf{f}_J\longmapsto \mathbf{m}_J,
\end{equation}
and similarly, its inverse $\MR^{-1}$.

Given a set of indices $\Lambda \subset \nabla^J$,
where $\nabla^J := \bigcup_{j=0}^J \nabla_j$ with 
$\nabla_0 := S_0$,
we define a {\it thresholding}
operator $\thr_{\Lambda}$
that leaves unchanged the components $d_{\lambda}$ 
of $\mathbf{m}_J$
if $\lambda \in \Lambda$,
and replaces it by $0$ otherwise.
Defining the level-dependent threshold values 
$(\epsilon_0,\epsilon_1,\ldots,\epsilon_J)$,
the set $\Lambda$ is given by 
\begin{equation}\label{eq3:MR_Lambda}
\lambda \in \Lambda \ {\rm if} \ 
\|d_{\lambda}\|_{L^p} \geq \epsilon_{|\lambda|}.
\end{equation}
Applying $\thr_{\Lambda}$ on the multi-scale decomposition 
$\mathbf{m}_J$ of $\mathbf{f}_J$ 
amounts then to building a multiresolution
approximation $\adap_{\Lambda}\mathbf{f}_J$ of $\mathbf{f}_J$, 
where the operator 
$\adap_{\Lambda}$ is given by
\begin{equation*}\label{eq3:adap_operator_def}
\adap_{\Lambda}:=\MR^{-1}\thr_{\Lambda}\MR,
\end{equation*}
in which all details of a certain level of regularity
have been discarded.

The multi-scale transform 
(\ref{eq3:M_U_M})
amounts to represent $\disc f_J$
in a wavelet space spanned by a biorthogonal wavelet basis.
Actually, as shown in \cite{Cohen03}, the cell-average
(\ref{eq3:average_finite_vol}) 
results from considering 
a scaling function  $\widetilde{\phi}_{\gamma}$ defined as 
\begin{equation}\label{eq3:MR_dual_scaling}
\widetilde{\phi}_{\gamma} := |\Omega_{\gamma}|^{-1} \chi_{\Omega_{\gamma}},
\end{equation}
where $\chi_{\Omega_{\gamma}}$ is a standard characteristic function
($\chi_{\Omega_{\gamma}}=1$ if $\x \in \Omega_{\gamma}$;
otherwise, $\chi_{\Omega_{\gamma}}=0$).
Therefore, the finite volume representation of $f(\x)$
on the grid $S_j$:
$\mathbf{f}_j=(f_{\gamma})_{\gamma \in S_j}$ 
can be equivalently defined with
$f_{\gamma} := \langle f,\widetilde{\phi}_{\gamma}\rangle$.
Similarly, 
introducing (\ref{eq3:MR_dual_scaling}) and (\ref{eq3:linear_prediction})
in (\ref{eq3:MR_detail}) defines a 
box wavelet $\widetilde{\psi}_\mu$
of order $\beta$:
\begin{equation}\label{eq3:wavelet_MR}
\widetilde{\psi}_{\mu} := \widetilde{\phi}_{\mu} - 
\sum_{\gamma \in R_I(\mu)}
\beta_{\mu,\gamma} \widetilde{\phi}_{\gamma}.
\end{equation}
Following (\ref{eq3:linear_prediction}) in this work 
we consider only average-interpolating wavelets
given by (\ref{eq3:wavelet_MR}), to generate dynamically
the adapted grids through multiresolution analysis.
Details (\ref{eq3:MR_detail})
can be defined as the coefficients 
related to $f$ when represented on a wavelet basis:
$d_{\mu} = \langle f,\widetilde{\psi}_\mu\rangle$.
Further details can be found in \cite{cohen2000a,muller2003}.

Based on \cite{Cohen03},
we can define the following $\ell^2$-norm:
\begin{equation*}\label{eq3:normalized_l2}
\|\mathbf{f}_J\|_2^2:= 2^{-dJ} \ds \sum_{\lambda \in S_J} (f_{\lambda})^2,
\end{equation*}
which corresponds to the $L^2$-norm of a piecewise constant function.
The following bound follows, 
\begin{equation}\label{eq3:adap_error_eps}
 \|\mathbf{f}_{J} - \adap_{\Lambda}\mathbf{f}_J \|_2 \leq C \eta_{\rm MR},
\end{equation}
as shown in Appendix~\ref{AppMR}
with the level-dependent threshold values:
\begin{equation}\label{eq3:epsilon_j}
\epsilon_j = 2^{d(j-J)/2}\eta_{\rm MR}, \quad j=0,1,\ldots,J,
\end{equation}
where $\eta_{\rm MR}$ corresponds to an accuracy tolerance.

\subsection{Poisson equation discretized on multiresolution grids}\label{subsec:PoissonMR}
Considering the following Poisson equation:
\begin{equation}\label{eq:Poisson}
\partial^2_\x\, V = f,
\end{equation}
with $\x\in \Omega$,
we can represent it on the finest finite volume grid $S_J$ as before
by taking cell-averages,
that is,
\begin{equation}\label{eq:f_J_Poisson}
\mathbf{f}_J=(\langle f,\widetilde{\phi}_{\gamma}\rangle)_{\gamma \in S_J}
= (\langle \partial^2_\x V,\widetilde{\phi}_{\gamma}\rangle)_{\gamma \in S_J}.
\end{equation}
Recall that $\mathbf{f}_J$ 
is an array of size $n=\#(S_J)$
(where $\#(\cdot)$ returns the cardinality of a set), 
$\mathbf{f}_J \in \R^n$,
corresponding to function $f$ discretized on the grid $S_J$.
Considering the space of square matrices of size $n$: $\M$,
we can define an operator
$\mathbf{A} \in \M$ such that
following (\ref{eq:f_J_Poisson}),
\begin{equation}\label{eq:AVf}
\mathbf{f}_J =\mathbf{A} \mathbf{V}_J + 
\Or\left( (\Delta x)^\alpha\right),
\end{equation}
where $\Delta x :={\rm diam}(\Omega_{\gamma}\vert_{\gamma \in S_J})$
corresponds to the spatial resolution of the finest grid $S_J$, and
$\mathbf{V}_J := (\langle V,\widetilde{\phi}_{\gamma}\rangle)_{\gamma \in S_J} \in \R^n$,
that is, the analytical solution $V$ to the Poisson equation
(\ref{eq:Poisson}) discretized on the grid $S_J$.
Operator $\mathbf{A}$ is no other than a spatial
discretization of the Laplace operator.
It is therefore a positive definite, and hence non-singular matrix
assuming appropriate boundary conditions at $\x \in \partial \Omega$
for the Poisson equation (\ref{eq:Poisson}).
In particular following (\ref{eq:AVf}),
the unique solution $\mathbf{V}_d \in \R^n$ of system
$\mathbf{A}\mathbf{V}_d=\mathbf{f}_J$
is an approximation of 
order $\alpha$ to $\mathbf{V}_J$.

Now, if we consider the multiresolution approximation
$\mathbf{f}^\epsilon_J:=\adap_{\Lambda}\mathbf{f}_J$
and $\mathbf{V}^\epsilon \in \R^n$, solution of the linear system:
$\mathbf{A}\mathbf{V}^\epsilon=\mathbf{f}^\epsilon_J$,
it can be shown that 
there is a constant $c>0$ such that the following bound holds:
\begin{equation}\label{eq:bound_analsol}
\|\mathbf{V}^\epsilon - \mathbf{V}_J\|_2 \leq c \left( (\Delta x)^\alpha
+ \eta_{\rm MR} \right).
\end{equation}
Given a finite volume spatial discretization of order $\alpha$,
the exact solution $\mathbf{V}_J$ of the
Poisson equation can be therefore approximated
according to a prescribed tolerance $\eta_{\rm MR}$,
even if the multiresolution analysis acts on the right-hand side
function.
In particular it follows that
the exact solution $\mathbf{V}_d$ of
the discrete Poisson equation 
$\mathbf{A}\mathbf{V}_d=\mathbf{f}_J$
is approximated by $\mathbf{V}^\epsilon$ in the same way
$\mathbf{f}^\epsilon_J$
does for 
$\mathbf{f}_J$:
\begin{equation*}\label{eq:hat_eps}
 \|\mathbf{V}_d - \mathbf{V}^\epsilon \|_2 \leq C\eta_{\rm MR}.
\end{equation*}

However,
the Laplacian will be discretized 
in practice
on an
adapted grid; 
therefore, we will not be 
solving system $\mathbf{A}\mathbf{V}^\epsilon=\mathbf{f}^\epsilon_J$.
Instead, an operator $\mathbf{\widetilde{A}}$ is introduced
which corresponds to the Laplacian
discretized on the adapted grid.
Denoting $\mathbf{\widetilde{V}} \in \R^n$, solution of system
$\mathbf{\widetilde{A}}\mathbf{\widetilde{V}}=\mathbf{f}^\epsilon_J$,
we can numerically demonstrate that bound (\ref{eq:bound_analsol}) 
remains valid, that is,
\begin{equation}\label{eq:res_coro}
\|\mathbf{\widetilde{V}} - \mathbf{V}_J\|_2 \leq c \left( (\Delta x)^\alpha
+ \eta_{\rm MR} \right),
\end{equation}
as long as $\mathbf{\widetilde{A}}$ is consistently defined 
within the multiresolution framework.
Therefore,
by applying the multiresolution analysis on
the right-hand side function and solving the discrete Poisson equation
on the corresponding adapted grid, we obtain a
solution $\mathbf{\widetilde{V}}$ that also verifies
\begin{equation}\label{eq:hat_tilde}
 \|\mathbf{V}_d - \mathbf{\widetilde{V}} \|_2 \leq C\eta_{\rm MR}.
\end{equation}

\section{Numerical implementation}\label{sec:imple}
We now describe 
the numerical technique conceived to construct
a Poisson solver
within the present multiresolution framework.
We consider the multiresolution finite volume implementation 
presented in \cite{Duarte11_SISC}.
For the sake of completeness some 
key aspects of this particular implementation
will be first recalled,
while more details and references
can be found in \cite{Duarte_Phd}.

\subsection{Construction of multiresolution grids}\label{subsec:MR}
The adapted grid is composed of a set
of nested dyadic grids: $S_j$, $j=0,1,\ldots,J$,
from the coarsest to the finest,
generated by refining recursively 
a given cell depending on the
local regularity of the time-dependent
variables,
measured by the details
at a given time.
Function $f$ in the Poisson equation
(\ref{eq:Poisson}) that
depends directly on these variables (and
hence varies also in time)
may be additionally considered if necessary
to generate the grids,
as well as the solution $V$
corresponding to the previous time-step.
These grids
are implemented in a 
multi-dimensional and Cartesian 
finite volume framework.
Data compression is achieved
by discarding the cells whose details
are not in $\Lambda$ according to 
(\ref{eq3:MR_Lambda}).
However, 
a {\it graded tree} $\Lambda_\epsilon$ is considered in practice
instead of $\Lambda$
because a certain data structure must be 
respected in order
to carry out the 
multiresolution transform $\MR$ in (\ref{eq3:M_U_M}).
In particular all cells in the interpolation stencils $R_I(\cdot)$
must be always available
(see \cite{Cohen03} for more details). 
Notice that $\Lambda \subset \Lambda_\epsilon$ and
error estimates like (\ref{eq3:adap_error_eps})
follows straightforwardly with $\adap_{\Lambda_\epsilon}$
instead of $\adap_{\Lambda}$.
Nevertheless, for the ease of reading we will keep the 
notation $\Lambda$ in the following to 
refer to a graded tree.

\begin{figure}[!ht]\label{fig:tree}
\centerline{\includegraphics[width=0.8\textwidth]{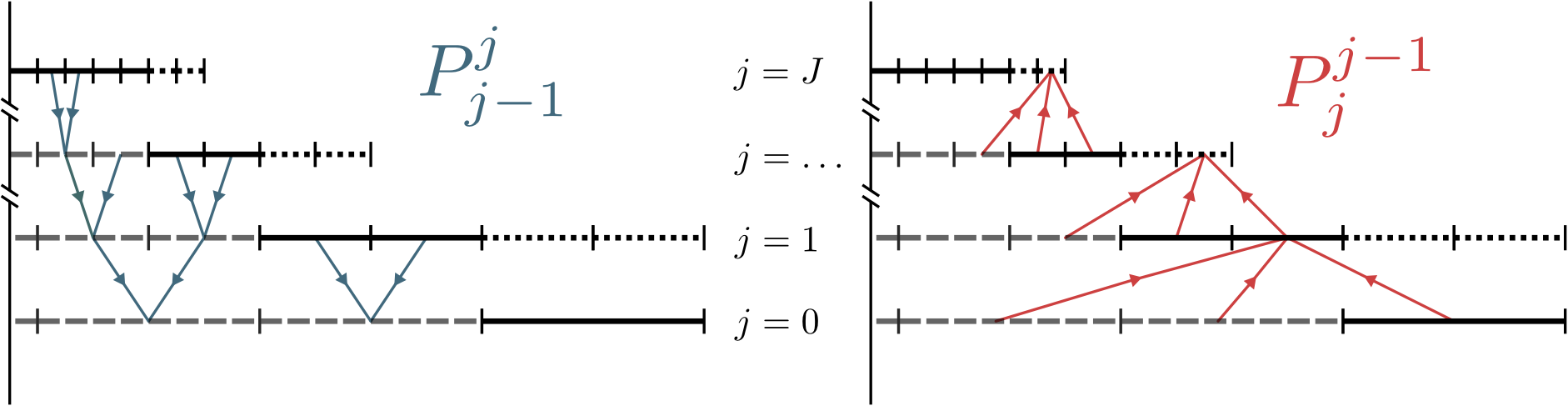}}
\caption{Part of a one-dimensional graded tree $\Lambda$, indicating
the leaves (solid lines) that form the adapted grid, as well as the
inner (dashed lines) and ghost (dotted lines) cells.
Projection $P^j_{j-1}$ (left) and prediction 
$P^{j-1}_j$ (right)
operations are also represented.}
\end{figure}

A graded tree-structure is used
to represent data in the computer memory
(see also \cite{Roussel03}).
Recalling the standard tree-structure terminology:
if $\Omega_{\mu} \subset \Omega_{\gamma}$ with
$|\mu| = |\gamma| + 1$, we say that $\Omega_\mu$ is a \textit{child} of 
$\Omega_\gamma$ and that $\Omega_\gamma$ is the \textit{parent} of $\Omega_\mu$. 
We thus define the \textit{leaves} 
$\leaf(\Lambda)$ of a 
\textit{tree} $\Lambda$ as the set of cells $\Omega_{\lambda}$,
$\lambda \in \leaf(\Lambda)$,
such that 
$\Omega_{\lambda}$ has no children in $\Lambda$.
The sets $\nabla_j$, $j=0,1,\ldots,J$, are distributed in 
$N_{\rm R}$ graded trees $\Lambda_r$, $r=1,\ldots,N_{\rm R}$,
where $N_{\rm R}:=N_{{\rm R}x}N_{{\rm R}y}N_{{\rm R}z}$, 
and
$N_{{\rm R}x}$, $N_{{\rm R}y}$, and
$N_{{\rm R}z}$ stand for the number of graded trees
or  {\it roots} per direction.
Denoting by $\tree(\Lambda _r)$ the set that 
contains the graded tree $\Lambda _r$ plus the
missing cells $\Omega_{\lambda}$ in the construction of sets
$\nabla_j$, we similarly
have that grid indices
$S_j$, $j=0,1,\ldots,J$, are distributed in 
$N_{\rm R}$ sets $\tree(\Lambda_r)$.
The adapted grid is thus given by sets
$\leaf(\Lambda_r)$, $r=1,\ldots,N_{\rm R}$,
with a total number of cells:
$N_\leaf= \sum_{r=1}^{N_{\rm R}}\#(\leaf(\Lambda_r))$.
If no adaptation is required: 
$\max N_{\rm L}=\#(S_J)=N_{{\rm R}}2^{dJ}$, 
that is, the size of the finest grid.
Ghost cells called {\it phantoms} are added 
to the adapted grid at level interfaces, in order
to always compute numerical fluxes at the
highest grid-level between two neighboring cells \cite{Roussel03}.
Cell-averages of phantoms are computed using the 
prediction operator (\ref{eq3:linear_prediction}); therefore,
the graded tree must also contain all cells needed 
to perform the corresponding interpolations.
Figure~\ref{fig:tree} depicts part of a one-dimensional graded tree
where the projection and prediction operators:
$P^j_{j-1}$ and $P^{j-1}_j$ according to 
(\ref{eq3:projection}) and
(\ref{eq3:linear_prediction}), respectively,
are schematically
described.

Following \cite{Cohen03}
a centered polynomial interpolation of accuracy order
$\beta = 2s+1$ is defined
for the projection operator (\ref{eq3:linear_prediction}), 
computed with the $s$ nearest neighboring cells
in each direction;
the procedure is exact for polynomials of degree $2s$.
In the numerical illustrations we will only consider the case
$\beta=3$ with one neighboring cell per direction ($s=1$)
including the diagonals in multidimensional configurations.
For the one-dimensional dyadic case (\ref{eq3:dyadic_1D}),
the latter is given by 
\begin{equation*}\label{eq3:polynomial_dyadic1D_order3}
\widehat{f}_{j+1,2k} = f_{j,k} + \frac{1}{8}(f_{j,k-1} - f_{j,k+1}), \qquad
\widehat{f}_{j+1,2k+1} = f_{j,k} + \frac{1}{8}(f_{j,k+1} - f_{j,k-1}),
\end{equation*}
as shown in Figure~\ref{fig:tree} (right).
In this case 
the set of coefficients 
$(\beta_{\mu,\gamma})_{\gamma \in R_I(\mu)}$ 
in (\ref{eq3:linear_prediction}) is given either by
$(1/8,1,-1/8)$ or by $(-1/8,1,1/8)$, regardless of the grid-level,
which are the scaling coefficients of a third-order
average-interpolating wavelet of type (\ref{eq3:wavelet_MR}).
Higher order formulae can be found in \cite{muller2003}, 
while extension to multi-dimensional 
Cartesian grids is easily obtained by a
tensorial product of the one-dimensional operator 
\cite{bihari1997,Roussel03}.
In general the interpolation stencil $R_I(\cdot)$ is given
by $(2s+1)^d$ cells.

Input parameters for the multiresolution implementation are: 
the maximum grid-level
$J$ corresponding to the finest spatial discretization;
the number of roots per direction
$N_{{\rm R}x}$, $N_{{\rm R}y}$, and
$N_{{\rm R}z}$;
and the threshold parameter
$\eta_{\rm MR}$ which defines the numerical
accuracy of the compressed representations following
(\ref{eq3:adap_error_eps}).

\subsection{Construction of the discrete Laplace operator}\label{subsec:const_lap}
Introducing the set
${\rm I}_\leaf := \{1, 2, \ldots, N_\leaf\}$,
we define a bijective function $h:D(h)\to{\rm I}_\leaf$,
with
\begin{equation*}
D(h):= \bigcup_{r=1}^{N_{\rm R}}  \leaf(\Lambda_r).
\end{equation*}
The set 
$\Theta_\leaf := (\Omega_{\lambda})_{h(\lambda)\in {\rm I}_\leaf }$
corresponds then to the adapted grid, defined by 
the leaves of the tree representation.
Multi-dimensional grids are thus arranged into 
a one-dimensional array $\Theta_\leaf$,
where each leaf is associated with a unique index from $1$ to 
$N_\leaf$ in ${\rm I}_\leaf$.

We then consider
for a given function 
$u(\x)$
and for each leaf
$\Omega_{\gamma} \in \Theta_\leaf$
($\gamma$ such that
$\gamma \in D(h)$)
the
following standard finite volume approximation:
\begin{equation}\label{eq:FV_flux}
\langle\partial^2_\x u,\widetilde{\phi}_{\gamma}\rangle =
|\Omega_{\gamma}|^{-1} \sum_{\mu \neq \gamma}
|\Gamma_{\gamma,\mu}| F_{\gamma,\mu}
+ \Or\left( \left[{\rm diam}(\Omega_{\gamma})\right]^\alpha\right),
\quad \gamma \in D(h),
\end{equation}
where $F_{\gamma,\mu}$ accounts for the flux across each interface
$\Gamma_{\gamma,\mu} := \overline{\Omega_{\gamma}} \cap 
\overline{\Omega_{\mu}}$.
Moreover, we can represent the flux computations by
\begin{equation}\label{eq:flux}
|\Omega_{\gamma}|^{-1}
\sum_{\mu \neq \gamma}
|\Gamma_{\gamma,\mu}| F_{\gamma,\mu}
=
\sum_{\lambda \in R_F(\gamma)} \alpha_{\gamma,\lambda} u_\lambda,
\end{equation}
where the {\it flux stencil}
$R_F(\gamma)$ is contained in one single grid-level 
($R_F(\gamma) \subset S_{|\gamma|}$)
and the set of coefficients $(\alpha_{\gamma,\lambda})_{\lambda \in R_F(\gamma)}$
establishes the order $\alpha$ of the approximation.
If the same scheme is considered throughout a given $S_j$,
then for any $\mu \neq \gamma$ such that $|\gamma|=|\mu|=j$ the set
of coefficients $(\alpha_{\gamma,\lambda})_{\lambda \in R_F(\gamma)}$ and 
$(\alpha_{\mu,\lambda})_{\lambda \in R_F(\mu)}$ are constant and component-wise
equal.
For instance, the classical centered second-order scheme in the 
one-dimensional dyadic case (\ref{eq3:dyadic_1D})
is given by
\begin{equation}\label{eq:1d_2ndorder_disc}
\langle\partial^2_x u,\widetilde{\phi}_{j,k}\rangle =
\Delta x_j^{-2} \left(
u_{j,k+1} - 2 u_{j,k} + u_{j,k-1}
\right)
+ 
\Or\left( \Delta x_j^2\right),
\end{equation}
where $\Delta x_j$ corresponds to the spatial resolution of grid $S_j$;
the set of coefficients $(\alpha_{\gamma,\lambda})_{\lambda \in R_F(\gamma)}$
in (\ref{eq:flux}) is thus given by
$\Delta x_{|\gamma|}^{-2}(1,-2,1)$.

The discrete Laplacian
$\mathbf{A} = (a_{i,l})_{i,l\in {\rm I}_\leaf}$ 
represented on the finest (uniform) finite volume grid $S_J$
is hence 
defined by setting
for each 
$i\in {\rm I}_\leaf$,
$\gamma = h^{-1}(i)$,
that is, for each 
leaf in $\Theta_\leaf$
($N_\leaf=\#(S_J)$):
\begin{equation}\label{eq:buildA}
a_{h(\gamma),h(\lambda)}
=
\alpha_{\gamma,\lambda},
\qquad
\forall
\lambda \in R_F(\gamma),
\end{equation}
and
\begin{equation}\label{eq:buildA_2_orig}
a_{h(\gamma),l} = 0,
\qquad
\forall l\in {\rm I}_\leaf \
{\rm s.t.} \
h^{-1}(l) \notin R_F(\gamma).
\end{equation}
In the case of 
(\ref{eq:1d_2ndorder_disc}),
the latter process 
(\ref{eq:buildA})--(\ref{eq:buildA_2_orig})
amounts to build the standard tridiagonal 
matrix with non-zero entries given by
$\Delta x_{|\gamma|}^{-2}(1,-2,1)$.
Nevertheless,
the finite volume flux representation (\ref{eq:FV_flux})
establishes that for a given interface 
$\Gamma_{\gamma,\mu}$ the following conservation property
is verified:
$F_{\gamma,\mu} +  F_{\mu,\gamma} = 0$.
Computing the flux $F_{\gamma,\mu}$ for $\Omega_\gamma$
amounts to evaluate also $F_{\mu,\gamma}$ for the neighboring
cell $\Omega_\mu$. 
Let us denote $F_{\gamma,\mu}^+$ as the right flux for 
$\Omega_\gamma$ and $F_{\mu,\gamma}^-$ as the left flux
for $\Omega_\mu$ along the normal direction to $\Gamma^+_{\gamma,\mu}$,
the right interface of $\Omega_\gamma$
(the same as the left interface of
$\Omega_\mu$: $\Gamma_{\mu,\gamma}^-$).
Similarly, $R_F^+(\gamma)$ stands for the stencil required to
compute $F_{\gamma,\mu}^+$ and naturally 
$R^-_F(\mu) \equiv R^+_F(\gamma)$.
Fluxes are then computed only once at each interface
and the same property is exploited to save computations
while constructing operator
$\mathbf{A}$.
Instead of (\ref{eq:buildA}), we can thus define
for each $i\in {\rm I}_\leaf$,
$\gamma = h^{-1}(i)$:
\begin{equation}\label{eq:buildA_1}
a_{h(\gamma),h(\lambda)}
=
a_{h(\gamma),h(\lambda)}+
\alpha_{\gamma,\lambda},
\qquad
\forall
\lambda \in R^+_F(\gamma),
\end{equation}
and
\begin{equation}\label{eq:buildA_3}
a_{h(\mu),h(\lambda)} = 
a_{h(\mu),h(\lambda)} 
-
\alpha_{\gamma,\lambda},
\quad
\forall \mu  \
{\rm s.t.} \
\Gamma^+_{\gamma,\mu} =
\overline{\Omega_{\gamma}} \cap 
\overline{\Omega_{\mu}};
\end{equation}
where initially all coefficients are set to zero, {\it i.e.},
$\mathbf{A} = 0$,
which automatically accounts for 
(\ref{eq:buildA_2_orig}).
For the example 
(\ref{eq:1d_2ndorder_disc}), we naturally obtain the same
tridiagonal matrix, but the coefficients
$(\alpha_{\gamma,\lambda})_{\lambda \in R^+_F(\gamma)}$
are now 
given by
$\Delta x_{|\gamma|}^{-2}(1,-1)$.
In general
the sparsity of the resulting matrix depends directly
on the stencil $R^+_F(\cdot)$ related to the flux computation scheme,
while the computational complexity of the procedure is
of $\Or(\#(S_J))$.

However,
we are interested in building the
Laplacian 
$\mathbf{\widetilde{A}} = (\widetilde{a}_{i,l})_{i,l\in {\rm I}_\leaf}$
represented on a multiresolution finite volume adapted grid,
meaning that $N_\leaf < \#(S_J)$.
The principle is the same, as we construct
$\mathbf{\widetilde{A}}$ by computing its elements following
(\ref{eq:buildA_1})--(\ref{eq:buildA_3}) 
with $\widetilde{a}_{i,l}$ instead of $a_{i,l}$.
Notice that for a given $\gamma$ such that
$|\gamma|=j$ all fluxes are computed at the
same grid $S_j$ in (\ref{eq:flux}).
In the case of adapted grids the latter involves
that fluxes are
computed on a locally uniform grid defined
by $R_F(\gamma)$.
Ghost cells are locally
introduced so that for a given  $\gamma$
all cells $\lambda \neq \gamma$ such that 
$\lambda \in R_F(\gamma)$ are available.
Given an adapted tree $\Lambda _r$,
let us denote by $\phan(\Lambda _r)$ the set of phantoms 
related to the tree $\Lambda _r$;
that is, 
all cells with index $\lambda$ such that for any
leaf $\Omega_\gamma$ in $\Theta_\leaf$,
$\lambda \in R_F(\gamma)$ but $\lambda \notin \tree(\Lambda _r)$.
Notice that by construction a phantom is always a child
of a leaf.
The variable values in these ghost cells are
computed based on the cells contained in the 
adapted representation $\tree(\Lambda _r)$,
as described in \S\ref{subsec:MR}.
Using the prediction operation (\ref{eq3:linear_prediction}),
variables at phantoms are defined by 
\begin{equation}\label{eq:phantom}
\widehat{u}_\mu = 
\sum_{\gamma \in R_I(\mu)}
\beta_{\mu,\gamma} u_\gamma,
\quad
|\mu|=|\gamma| +1,
\end{equation}
such that
\begin{equation*}\label{eq:proj_ghost}
  u_{\gamma} = |\Omega_{\gamma}|^{-1} 
\sum_{|\mu|=|\gamma| +1,\Omega_{\mu} 
\subset \Omega_{\gamma}} |\Omega_{\mu}|\widehat{u}_{\mu},
\quad \gamma \in \leaf(\Lambda_r).
 \end{equation*}
Recalling that 
a phantom stands at the place of a discarded cell, 
we have that 
$\widehat{u}_\mu$ involves an approximation
error of $\Or(\epsilon_{|\mu|})$ 
according to (\ref{eq3:MR_detail}),
and the multiresolution
error framework remains perfectly valid.
Moreover, this construction guarantees a consistent and conservative
representation at inter-grid interfaces.
\begin{figure}[!ht]\label{fig:tree_operations}
\centerline{(a)
\hspace{0.4\textwidth}\hfill
(b)\hspace{0.4\textwidth}}
\centerline{
\includegraphics[width=0.4\textwidth]{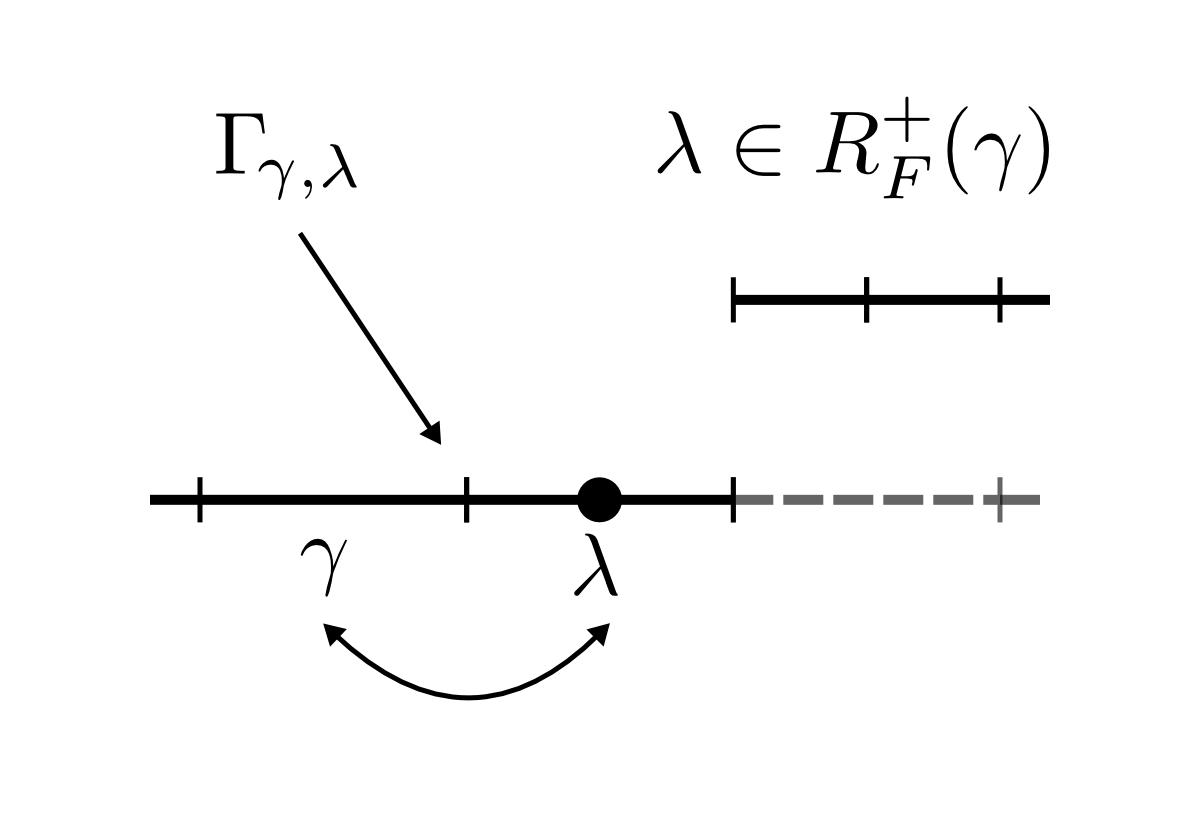}\hfill
\includegraphics[width=0.4\textwidth]{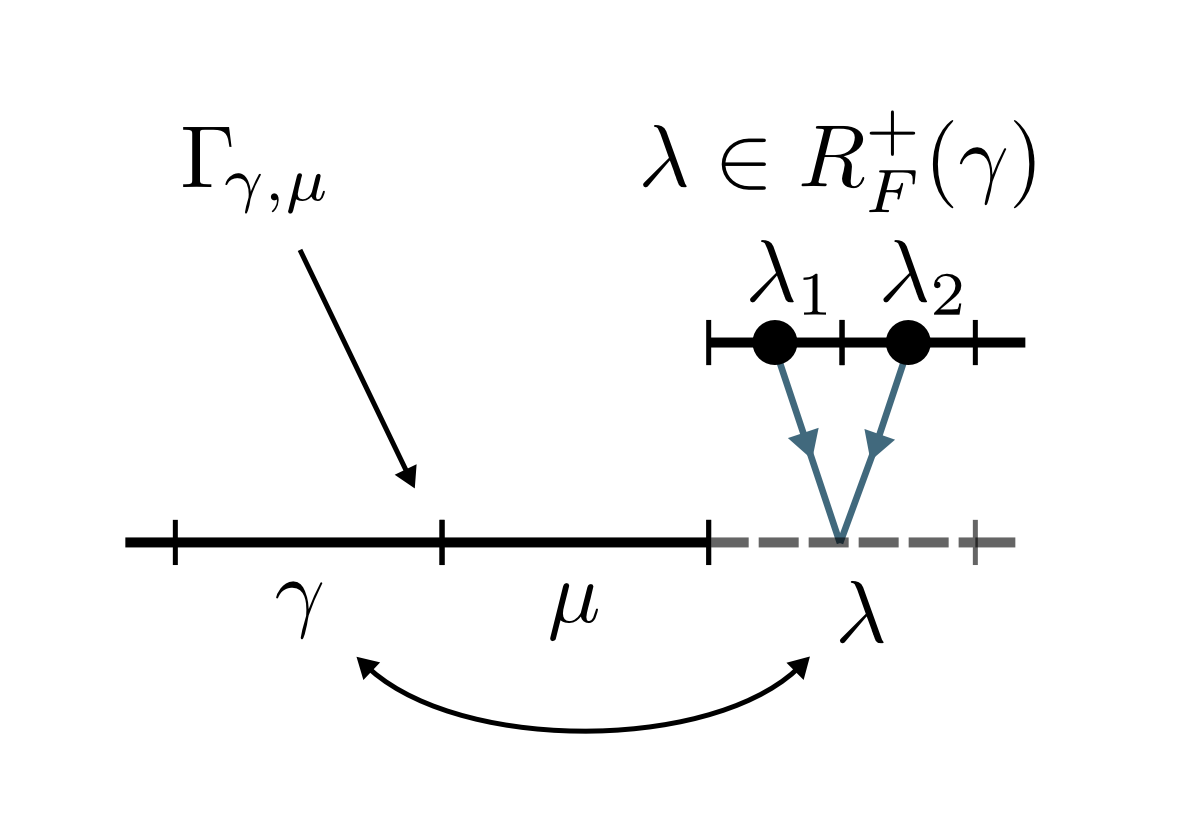}}
\centerline{(c)
\hspace{0.4\textwidth}\hfill
(d)\hspace{0.4\textwidth}}
\centerline{
\includegraphics[width=0.4\textwidth]{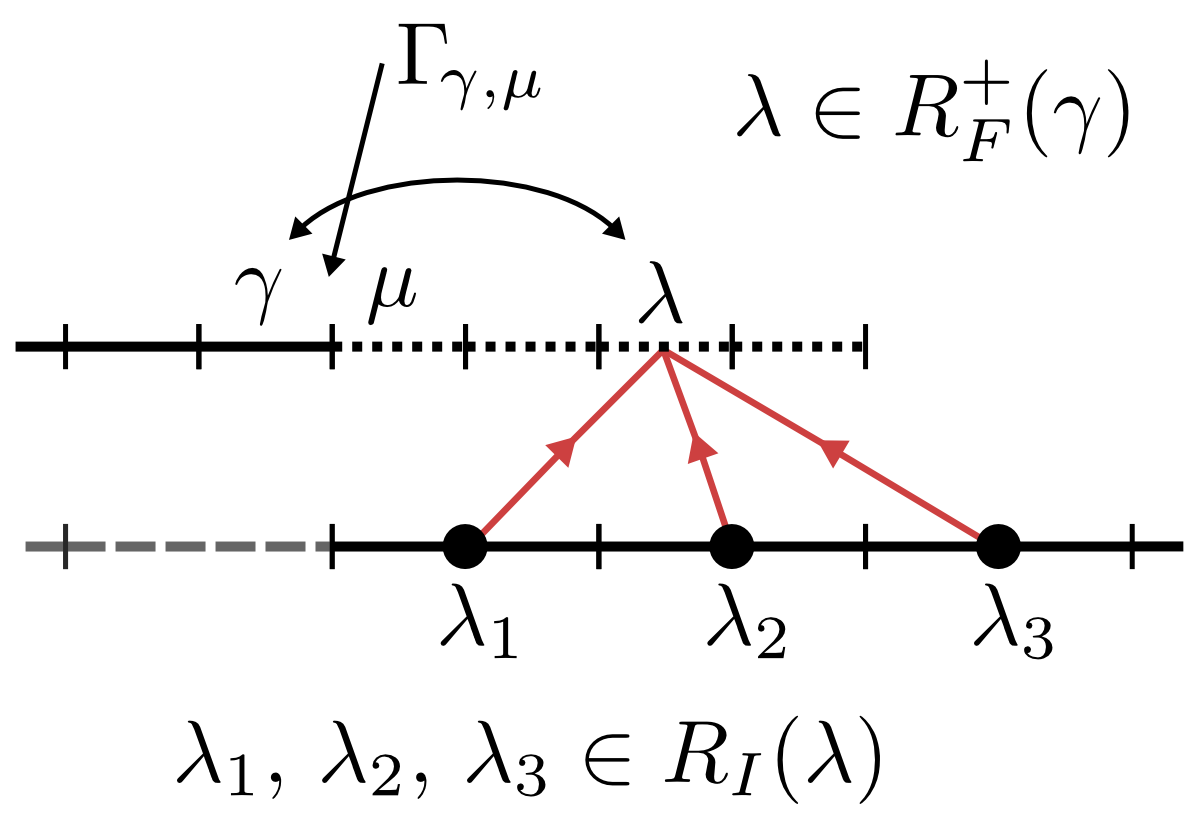}\hfill
\includegraphics[width=0.4\textwidth]{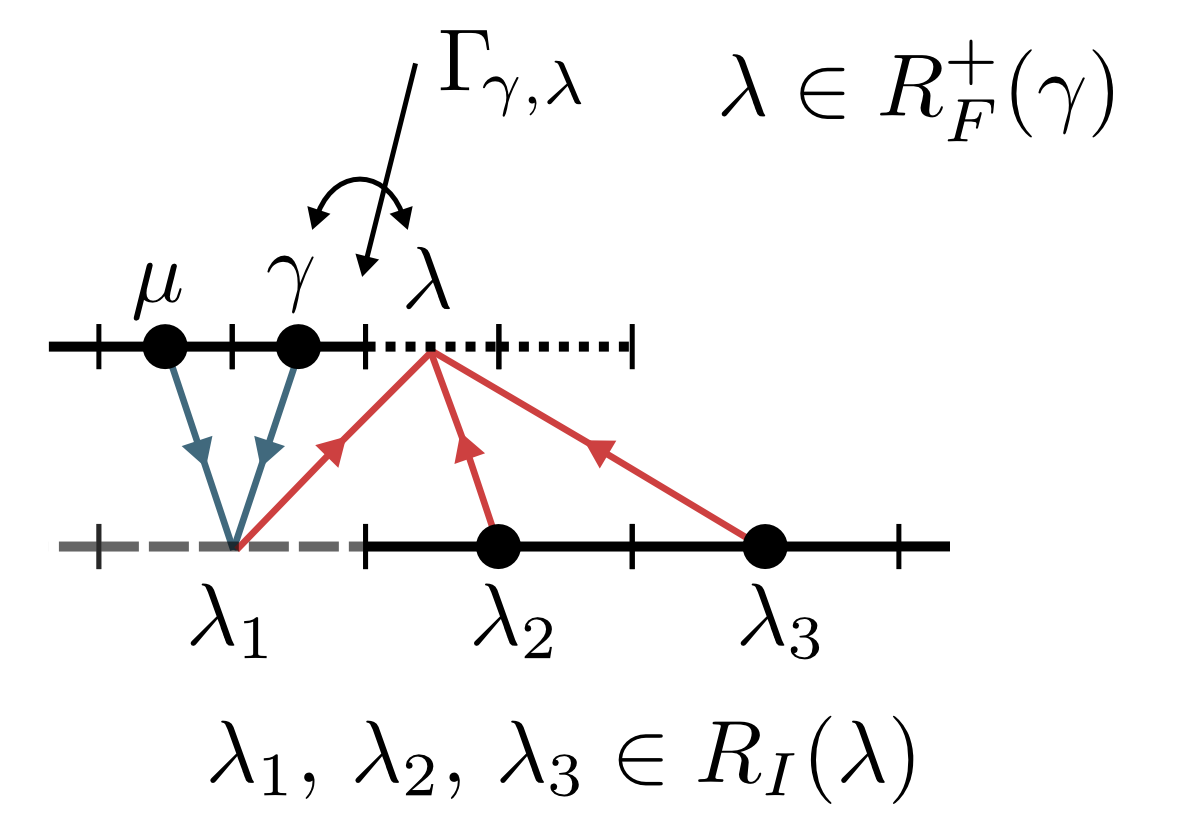}}
\caption{Computation of coefficients at inter-grid interfaces
when $\Omega_\lambda$ is contained in the flux stencil $R^+_F(\gamma)$
of a leaf $\Omega_\gamma$
and $\Omega_\lambda$ is a leaf (a), an inner cell (b), or 
a phantom (c)--(d).
Leaves, inner cells, and phantoms are represented with solid, dashed,
and dotted lines, respectively.
Coefficients related to cell $\Omega_\lambda$ are written
in terms of leaves marked with $\bullet$.}
\end{figure}

\begin{figure}[!ht]\label{fig:tree_operations2}
\centerline{(a)
\hspace{0.4\textwidth}\hfill
(b)\hspace{0.4\textwidth}}
\centerline{
\includegraphics[width=0.4\textwidth]{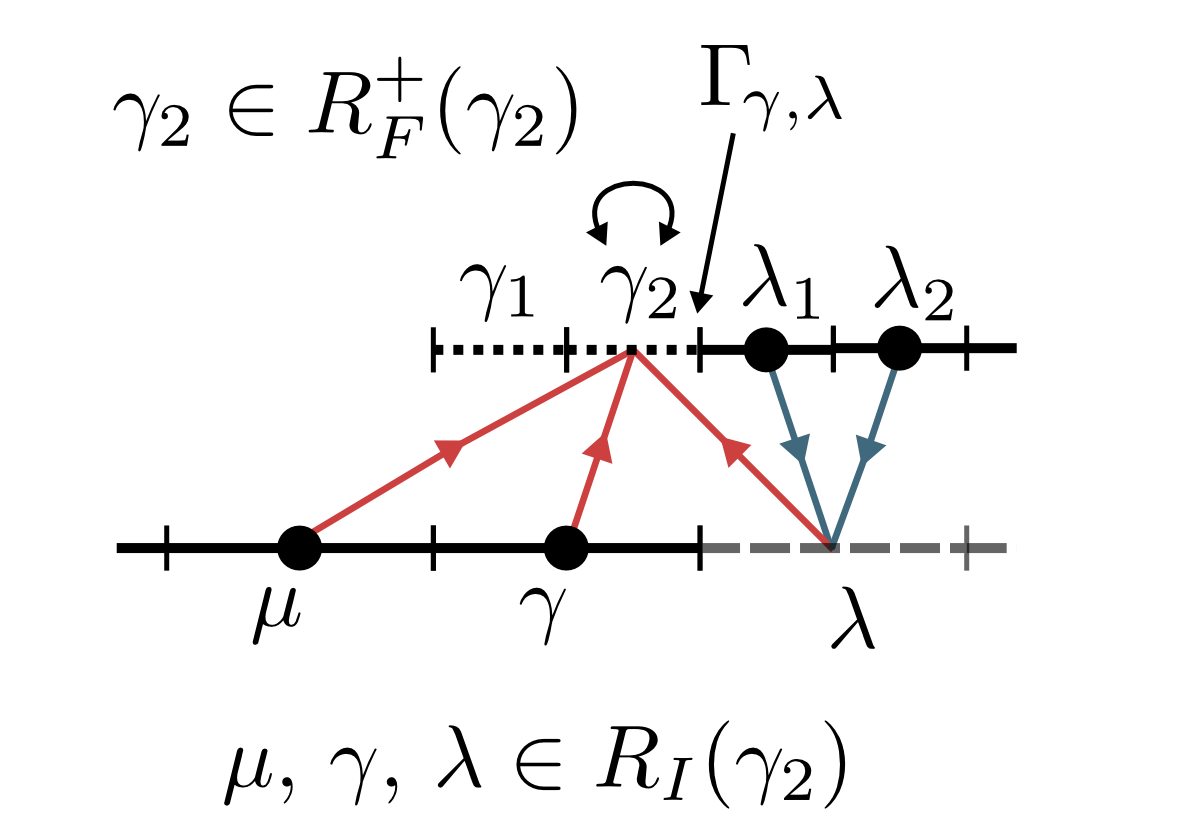}\hfill
\includegraphics[width=0.4\textwidth]{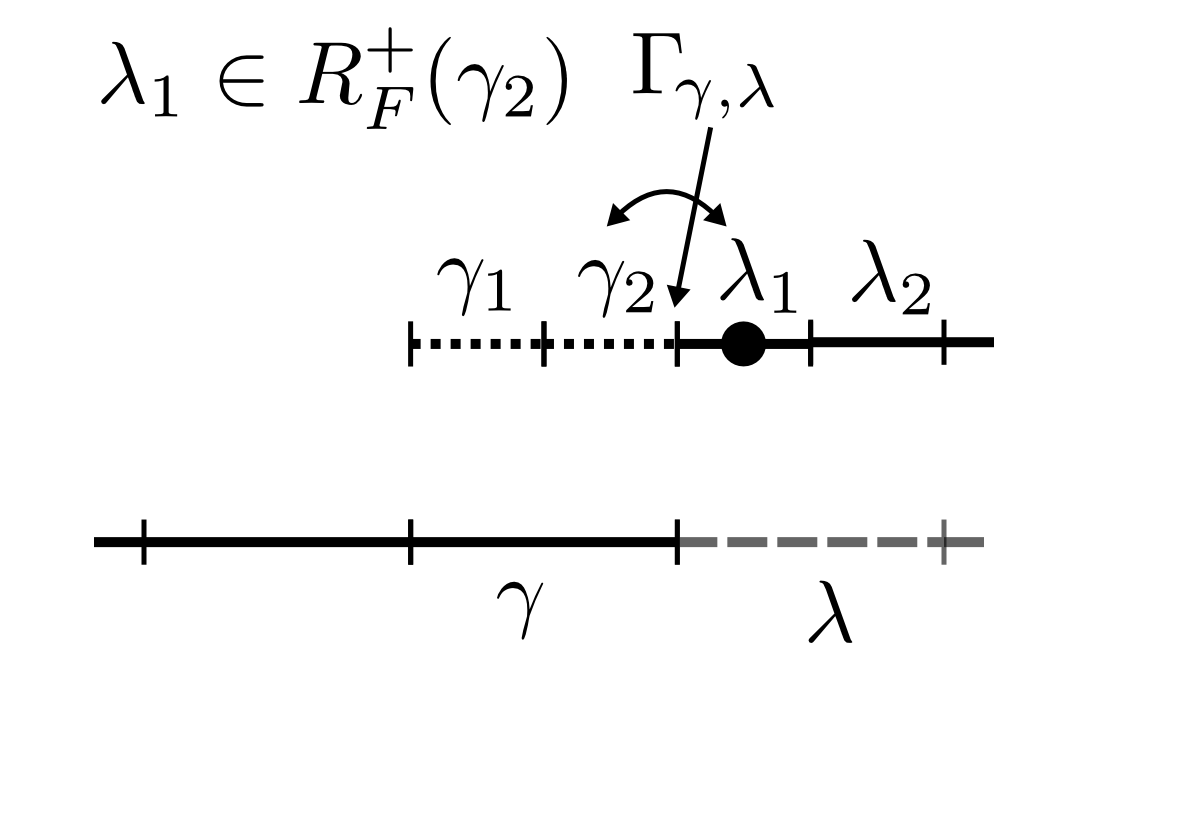}}
\caption{
Computation of coefficients at inter-grid interfaces
when a leaf $\Omega_\gamma$ shares an interface with another one of
higher resolution.
Leaves, inner cells, and phantoms are represented with solid, dashed,
and dotted lines, respectively.
Coefficients related to cell 
$\Omega_\gamma$ (a) and
$\Omega_{\lambda_1}$ (b) are written
in terms of leaves marked with $\bullet$.}
\end{figure}

Given a certain leaf $\Omega_\gamma$, for 
each cell $\Omega_\lambda$ within the flux stencil $R^+_F(\gamma)$
($\lambda \in R^+_F(\gamma)$) there are three cases:
\begin{enumerate}[(i)]
 \item \label{eq:leaf_a} Cell $\Omega_\lambda$ is a leaf and thus belongs to the adapted grid, {\it i.e.},
$\lambda \in \bigcup_{r=1}^{N_{\rm R}}  \leaf(\Lambda_r)$.
The coefficient $a_{h(\gamma),h(\lambda)}$ is computed according to 
(\ref{eq:buildA_1}) as for a uniform grid (see, for instance, Figure \ref{fig:tree_operations}(a)).
\item \label{eq:tree_a}
Cell $\Omega_\lambda$ belongs to the set of adapted grids but it is not
 a leaf, {\it i.e.},
 $\lambda \in \bigcup_{r=1}^{N_{\rm R}}  \tree(\Lambda_r) \wedge
 \lambda \notin \bigcup_{r=1}^{N_{\rm R}}  \leaf(\Lambda_r)$.
In this case inner cells are linked to leaves
using the projection operator (\ref{eq3:projection}).
In Figure \ref{fig:tree_operations}(b),
(\ref{eq:buildA_1}) is replaced by
\begin{align*}
&a_{h(\gamma),h(\lambda_i)}
=
a_{h(\gamma),h(\lambda_i)}+
|\Omega_{\lambda}|^{-1} 
|\Omega_{\lambda_i}|
\alpha_{\gamma,\lambda},
\qquad
\Omega_{\lambda_i} 
\subset \Omega_{\lambda},
i=1,2.
\end{align*}
 \item \label{eq:phantom_a} Cell $\Omega_\lambda$ is a phantom, {\it i.e.},
 $\lambda \in \bigcup_{r=1}^{N_{\rm R}}  \phan(\Lambda_r)$.
As established by (\ref{eq:phantom}),
phantoms are linked to leaves 
and/or inner cells using the 
prediction operation (\ref{eq3:linear_prediction}).
For Figure \ref{fig:tree_operations}(c)
we thus have
\begin{align*}
&a_{h(\gamma),h(\lambda_i)}
=
a_{h(\gamma),h(\lambda_i)}+
\beta_{\lambda,\lambda_i}
\alpha_{\gamma,\lambda},
\qquad
\lambda_i \in R_I(\lambda),
i=1,2,3;
\end{align*}
whereas in Figure \ref{fig:tree_operations}(d)
$\Omega_{\lambda_1}$ is not a leaf and hence,
$a_{h(\gamma),h(\lambda_1)}$ must be replaced by
\begin{align*}
&a_{h(\gamma),h(\mu)}
=
a_{h(\gamma),h(\mu)}+
|\Omega_{\lambda_1}|^{-1} 
|\Omega_{\mu}|
\beta_{\lambda,\lambda_1}
\alpha_{\gamma,\lambda},
\qquad
\Omega_{\mu} 
\subset \Omega_{\lambda_1},
\lambda_1 \in R_I(\lambda),
\\
&a_{h(\gamma),h(\gamma)}
=
a_{h(\gamma),h(\gamma)}+
|\Omega_{\lambda_1}|^{-1} 
|\Omega_{\gamma}|
\beta_{\lambda,\lambda_1}
\alpha_{\gamma,\lambda},
\qquad
\Omega_{\gamma} 
\subset \Omega_{\lambda_1},
\lambda_1 \in R_I(\lambda),
\end{align*}
combining both inter-level operations
(\ref{eq3:projection}) and
(\ref{eq3:linear_prediction}).
\end{enumerate}
Additionally, 
if a leaf $\Omega_\gamma$ shares
an interface $\Gamma_{\gamma,\lambda}$
with another of higher resolution,
coefficients $\widetilde{a}_{h(\gamma),l}$ 
are computed at grid-level $|\gamma|+1$,
considering the corresponding phantoms,
children of $\Omega_\gamma$, at $\Gamma_{\gamma,\lambda}$.
For instance,
when considering the interface $\Gamma_{\gamma,\lambda}$
in Figure \ref{fig:tree_operations2},
the stencil 
$R^+_F(\gamma_2)$ is considered corresponding to 
the phantom $\Omega_{\gamma_2}$, child of $\Omega_{\gamma}$.
For $\Omega_{\gamma_2}$ contained in $R^+_F(\gamma_2)$ 
(Figure \ref{fig:tree_operations2}(a)),
we thus have
\begin{align*}
a_{h(\gamma),h(\mu)}
=
a_{h(\gamma),h(\mu)}+
|\Omega_{\gamma}|^{-1} 
|\Omega_{\gamma_2}|
\beta_{\gamma_2,\mu}
\alpha_{\gamma_2,\gamma_2},
\qquad &
\Omega_{\gamma_2} 
\subset \Omega_{\gamma},
\mu \in R_I(\gamma_2),
\\
a_{h(\gamma),h(\gamma)}
=
a_{h(\gamma),h(\gamma)}+
|\Omega_{\gamma}|^{-1} 
|\Omega_{\gamma_2}|
\beta_{\gamma_2,\gamma}
\alpha_{\gamma_2,\gamma_2},
\qquad &
\Omega_{\gamma_2} 
\subset \Omega_{\gamma},
\gamma \in R_I(\gamma_2),
\\
a_{h(\gamma),h(\lambda_i)}
=
a_{h(\gamma),h(\lambda_i)}+
|\Omega_{\gamma}|^{-1} 
|\Omega_{\gamma_2}|
\beta_{\gamma_2,\lambda}
|\Omega_{\lambda}|^{-1} 
|\Omega_{\lambda_i}|
\alpha_{\gamma_2,\gamma_2},
\qquad &
\Omega_{\gamma_2} 
\subset \Omega_{\gamma},
\lambda \in R_I(\gamma_2),
\Omega_{\lambda_i} 
\subset \Omega_{\lambda},
i=1,2;
\end{align*}
whereas for 
$\Omega_{\lambda_1}$, also contained in $R^+_F(\gamma_2)$ 
(Figure \ref{fig:tree_operations2}(b)):
\begin{align*}
& a_{h(\gamma),h(\lambda_1)}
=
a_{h(\gamma),h(\lambda_1)}+
|\Omega_{\gamma}|^{-1} 
|\Omega_{\gamma_2}|
\alpha_{\gamma_2,\lambda_1},
\qquad 
\Omega_{\gamma_2} 
\subset \Omega_{\gamma}.
\end{align*}

The general, multi-dimensional
scheme 
to 
construct the discrete Laplacian
is detailed
for the sake of completeness in 
Appendix~\ref{AppScheme}.
The algorithm (with computational complexity 
$\Or(N_{\rm L})$)
considers multiplications
and combinations of constant coefficients
coming from the finite volume discretization
and the inter-level multiresolution operations
as previously shown.
The operator
$\mathbf{\widetilde{A}}$ is thus represented 
directly on the adapted grid strongly 
coupling consistent inter-grid and conservation 
properties.
In this study we do not develop our own linear
solver to solve the discrete Poisson equation.
We rather rely on dedicated solvers available in the literature
whose performance to solve multidimensional linear systems
have already been demonstrated.
Consequently, 
matrix $\mathbf{\widetilde{A}}$
must be an input to these solvers;
otherwise, the same operations 
described here
would have been performed
without having to save the matrix entries in memory.
Here operator $\mathbf{\widetilde{A}}$
is stored using a standard CSR (Compressed
Sparse Row) format for sparse matrices.

\section{Streamer discharge simulations}\label{sec:streamer}
Classical fluid model for streamers in air 
at atmospheric pressure 
is given by drift-diffusion equations
consistently coupled with a Poisson equation \cite{Babaeva:1996,Kulikovsky:1997c}:
\def\SeP{\ne\nu_{\rm i}}
\def\SpP{\ne\nu_{\rm i}}
\def\SeM{\ne(\nu_{{\rm a}2}+\nu_{{\rm a}3}) -  \ne\np\beta_{\rm ep}}
\def\SpM{\ne\np\beta_{\rm ep} - \nn\np\beta_{\rm np}}
\def\SnM{\nn\np\beta_{\rm np}}
\def\SnP{\ne(\nu_{{\rm a}2}+\nu_{{\rm a}3})}
\def\SD{\nn\gamma}
\begin{equation}\label{trasp} 
\left.
\begin{array}{l}
\partial_t \ne +\Dx \cdot(\ne\,\ve) 
  -\Dx\cdot(\De\ \Dx\ne) =     
  \SeP  -\SeM + \SD + \Sph, \\[1.5ex]
\partial_t \np +\Dx\cdot(\np\vp) 
  -\Dx\cdot(\Dp\,\Dx\np) = 
  \SpP - \SpM + \Sph, \\[1.5ex]
\partial_t \nn  +\Dx\cdot(\nn\vn) 
  -\Dx\cdot(\Dn\,\Dx\nn)  =    
  \SnP -\SnM - \SD, 
\end{array} 
\right\}
\end{equation}
\begin{equation}
\varepsilon_0\, \Dx \cdot \vec E  = -q_{\rm e}(\np-\nn-\ne), \quad \vec E = - \Dx \phi,  \label{poisson}
\end{equation}
where $\x\in \R^d$,
$n_{\rm i}$ is the  density of charged species ${\rm i}$ 
(e: electrons, p: positive ions, n: negative ions),
$\phi$  and $\vec E$
stand, respectively, for the electric potential and field,
and $\vec{v}_{\rm i}= \mu_{\rm i} \vec E$ is the drift velocity.
We denote by
$D_{\rm i}$ and $\mu_{\rm i}$, respectively, the
diffusion coefficient and the mobility of charged species ${\rm i}$,
$q_{\rm e}$ is the absolute value of the electron charge,
and   $\varepsilon_0$  is the permittivity of free space.
Moreover,
$\nu_{\rm i}$ is the electron impact ionization coefficient, $\nu_{{\rm a}2}$ and $\nu_{{\rm a}3}$
are the two-body and three-body electron attachment coefficients, $\beta_{\rm ep}$
and $\beta_{\rm np}$
are, respectively, the electron-positive ion
and
negative-positive ion recombination coefficients,
and $\gamma$ is the detachment coefficient.
All these coefficients depend
on the local reduced electric field $E/N_{\rm air}$
and thus vary in time and space,
where $E=|\vec E|$ is the electric field magnitude and $N_{\rm air} $ is the
air neutral density.
For test studies presented in this paper, the transport parameters for
air are taken from \cite{Morrow:1997}; detachment and attachment coefficients, 
respectively, from \cite{Benilov:2003,Kossyi:1992}; and other reaction rates, also
from \cite{Morrow:1997}. Diffusion coefficients for ions are derived from mobilities using classical Einstein relations.
Our reference density for air is $N_{\rm air}=2.688\times 10^{19}\,$cm$^{-3}$.

For positive streamers
a sufficient number of seed-electrons needs to be present in front of the streamer head
as the direction of electron motion is opposed to the streamer propagation 
(see \cite{Bourdon:2007} and references therein). 
Photoionization is in general an accepted mechanism to produce such seed-electrons
in nitrogen-oxygen mixtures.
It is therefore introduced into the
drift-diffusion system (\ref{trasp}) as a source term ($\Sph$)
that needs to be evaluated in general at each time-step
for all points of the computational domain.
Computation of $\Sph$ is detailed in 
Appendix~\ref{PhotoIonSec}
which requires the iterative solution of six elliptic equations given by
(\ref{SP3}) with boundary conditions (\ref{SP3BC}).
Iterating three times amounts then to solve 18 elliptic equations
per time-step.

In what follows we will first assess the theoretical
validity of the mathematical description conducted in \S\ref{subsec:PoissonMR}
and the numerical implementation described in \S\ref{subsec:const_lap}.
This study will be conducted on a simplified multi-dimensional model
with known analytical solution that mimics the spatial configuration
typically found in streamer discharges.
In a second part 
we will present
two-dimensional
double-headed streamer simulations
modeled by (\ref{trasp})--(\ref{poisson}),
for which we will evaluate the performance of
different linear solvers implemented to solve the discrete Poisson equations.
Finally,
dynamic grid adaptation will be analyzed for
the numerical simulation of
two interacting positive streamers
in a two-dimensional configuration that leads to streamer merging.

\subsection{Numerical validation}
We first investigate the validity of 
bound (\ref{eq:res_coro}) (and (\ref{eq:hat_tilde})).
That is, 
the numerical error related to grid adaptation and
data compression is of $\Or(\eta_{\rm MR})$,
where $\eta_{\rm MR}$ is a user-defined accuracy tolerance.
Given a set of constant parameters:
$a$, $b$, and $\sigma$, let us consider an exponential function $\phi(\x)$
on a multi-dimensional domain 
$\Omega \subset \R^d$,
\begin{equation}\label{SOLUTION}
 \phi(\x) = g(\x) + b =  a \exp \left(-|\x|^2/\sigma^2 \right) + b, \quad
 \x \in \Omega
\end{equation}
that verifies the following Poisson equation:
\begin{equation}\label{PTestCase1}
\Dx^2 \phi(\x) =  \rho(\x), \quad
\rho(\x)  = \frac{4}{\sigma^2}\left(\frac{ |\x|^2}{\sigma^2} - 1\right) g(\x), \quad
\x \in \Omega,
\end{equation}
with boundary conditions,
\begin{equation}\label{BC}
\phi(\x) =g(\x)+b, \quad \x \in  \partial\Omega.
\end{equation}
Using the standard, second-order centered scheme 
(similar to (\ref{eq:1d_2ndorder_disc})),
we discretize equation (\ref{PTestCase1})
on a two-
and a three-dimensional
region:  $[-0.5, 0.5]^2$ and
$[-0.5, 0.5]^3$, respectively, and we
consider the set of parameters:
$a=10$, $b=20$, and $\sigma=0.005$. 
The value of $\sigma$ 
has been chosen such that function
$\rho(\x)$ 
exhibits similar steep gradients as those found in
a developed streamer head
modeled by (\ref{trasp})--(\ref{poisson}).
Since $g(\x)$ decays rapidly toward the boundaries,
we consider
Dirichlet boundary conditions in (\ref{BC}):
$\phi(\x)=b$, whereas symmetric boundary conditions
are taken in order to consider a reduced
two- and a three-dimensional computational
domain: $[-0.5, 0.5] \times [0, 0.5]$ and
$[0, 0.5]^3$, respectively.
\begin{figure}[!ht] 
\centering
\includegraphics[width=0.49\hsize]{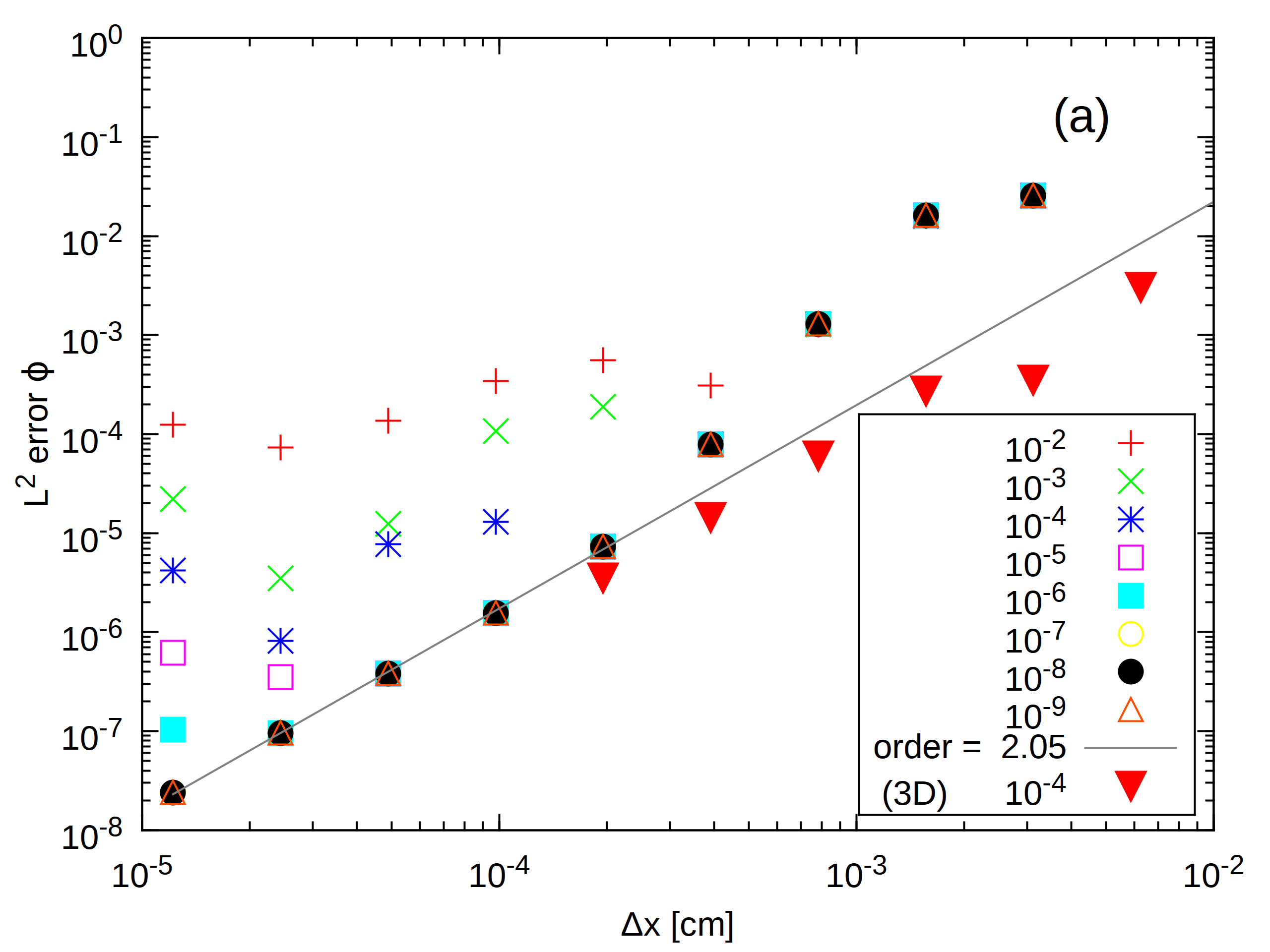}
\includegraphics[width=0.49\hsize]{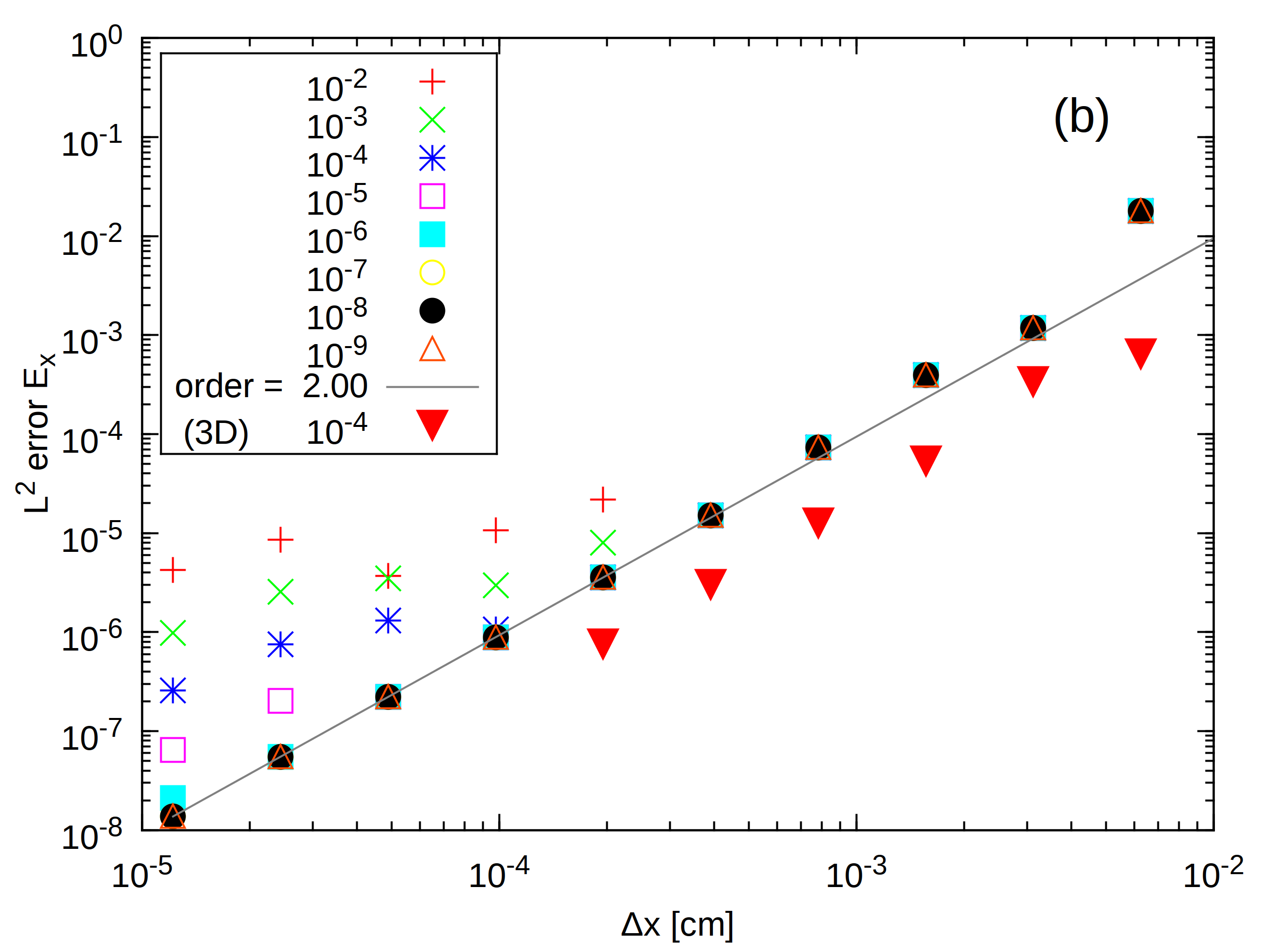}
\caption{$L^2$-errors between analytical and numerical solutions
$\phi$ of Poisson equation (\ref{PTestCase1}) (a)
and component $E_x$ of $\vec E = -\Dx \phi$ (b) 
for several threshold values $\eta_{\rm MR}$
( $\eta_{\rm MR}=10^{-4}$ for the three-dimensional problem).}
\label{order-phi} 
\end{figure}  

Figure \ref{order-phi} shows normalized $L^2$-errors 
between the analytical solution (\ref{SOLUTION}) 
and the numerical solution of the Poisson equation 
(\ref{PTestCase1})
discretized 
on an adapted grid, obtained with several threshold values $\eta_{\rm MR}$.
The resulting linear systems were solved using MUMPS\footnote{
Release 4.10.0. 
MUMPS home page: {\tt http://graal.ens-lyon.fr/MUMPS/}}
\cite{mumps,mumps1},
a direct linear system solver,
and BoomerAMG \cite{BoomerAMG} 
(contained in the \hypre library\footnote{
Release 2.8.0b. 
\hypre home page: 
{\tt http://www.llnl.gov/CASC/linear\_solvers/}}),
an iterative solver, for the two- and 
three-dimensional problems, respectively.
The finest spatial discretization is denoted by $\Delta x$,
and it is set by the choice of the maximum level $J$
in the multiresolution analysis
and the number of roots per direction:
$N_{{\rm R}x}$, $N_{{\rm R}y}$,  and $N_{{\rm R}z}$.
For the two-dimensional case, $\Delta x = 1/(N_{{\rm R}x}2^J)$ with $J=5,6,\ldots,13$,
$N_{{\rm R}x}=10$, and $N_{{\rm R}y}=5$;
whereas for the three-dimensional one:
$\Delta x = 0.5/(N_{{\rm R}x}2^J)$ with $J=4,5,\ldots,9$,
$N_{{\rm R}x}=N_{{\rm R}y}=N_{{\rm R}z}=5$.
All computations were performed on a work station
with 24 GB of computer memory.
For streamer discharge simulations an accurate resolution of the electric field: 
$\vec{E} = -\Dx \phi$, is essential for good physical descriptions. 
Therefore, we have also computed $\vec E$ with a
second-order, centered approximation,
and compared it against its analytical counterpart:
$\vec{E}= 2\x\, g(\x)/\sigma^2$.
In both cases, for $\phi$ and $\vec E$,
the numerical errors behave like a second
order spatial approximation even if the solutions are
computed on an adapted grid,
especially for relatively coarse discretizations
or sufficiently
fine multiresolution threshold values.
For finer resolutions, the numerical errors coming 
from the adaptive multiresolution become more
dominant and
the numerical errors are
effectively bounded by the threshold
parameter $\eta_{\rm MR}$.
Bounds (\ref{eq:res_coro}) and (\ref{eq:hat_tilde})
prove then to describe accurately the behavior
of the numerical approximations
when solving a Poisson equation on a multiresolution
adapted grid.
\begin{figure}[!ht] 
\centering
\includegraphics[width=0.39\hsize]{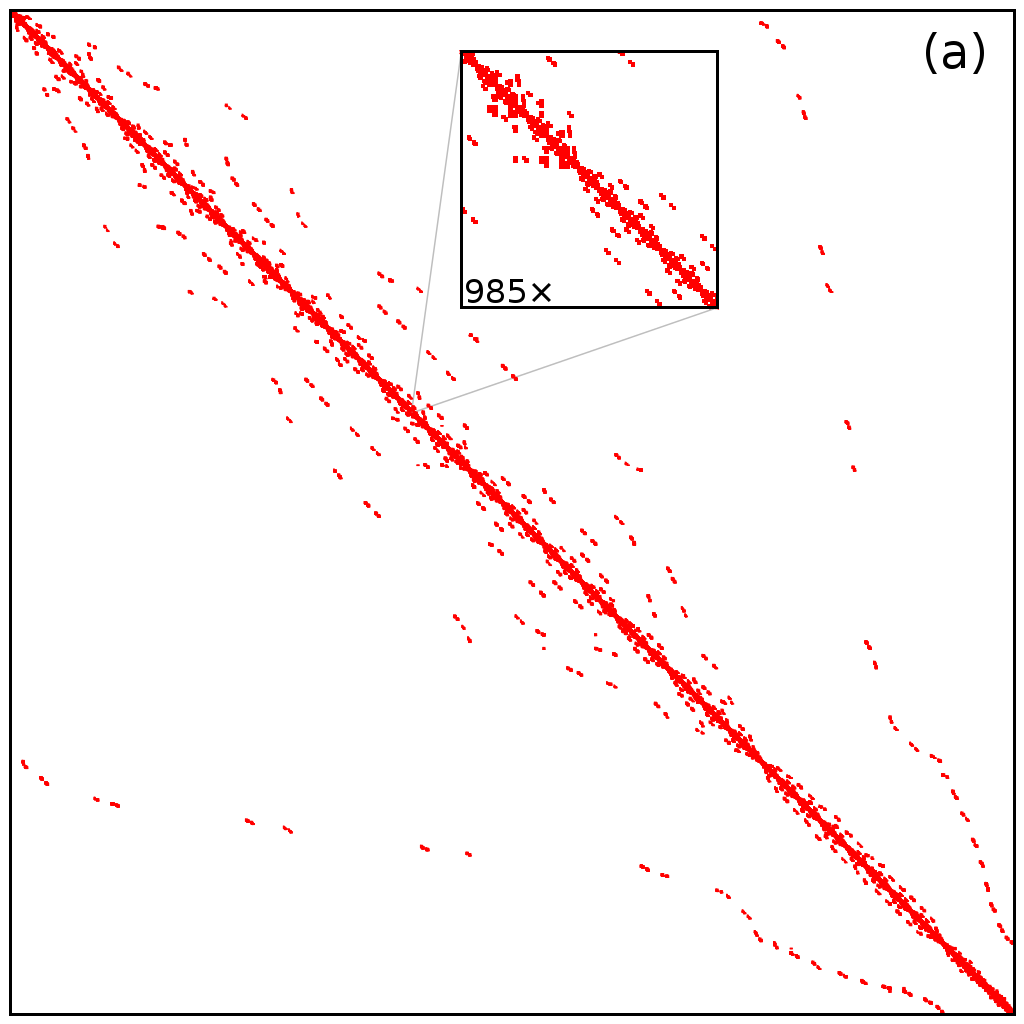}
\includegraphics[width=0.49\hsize]{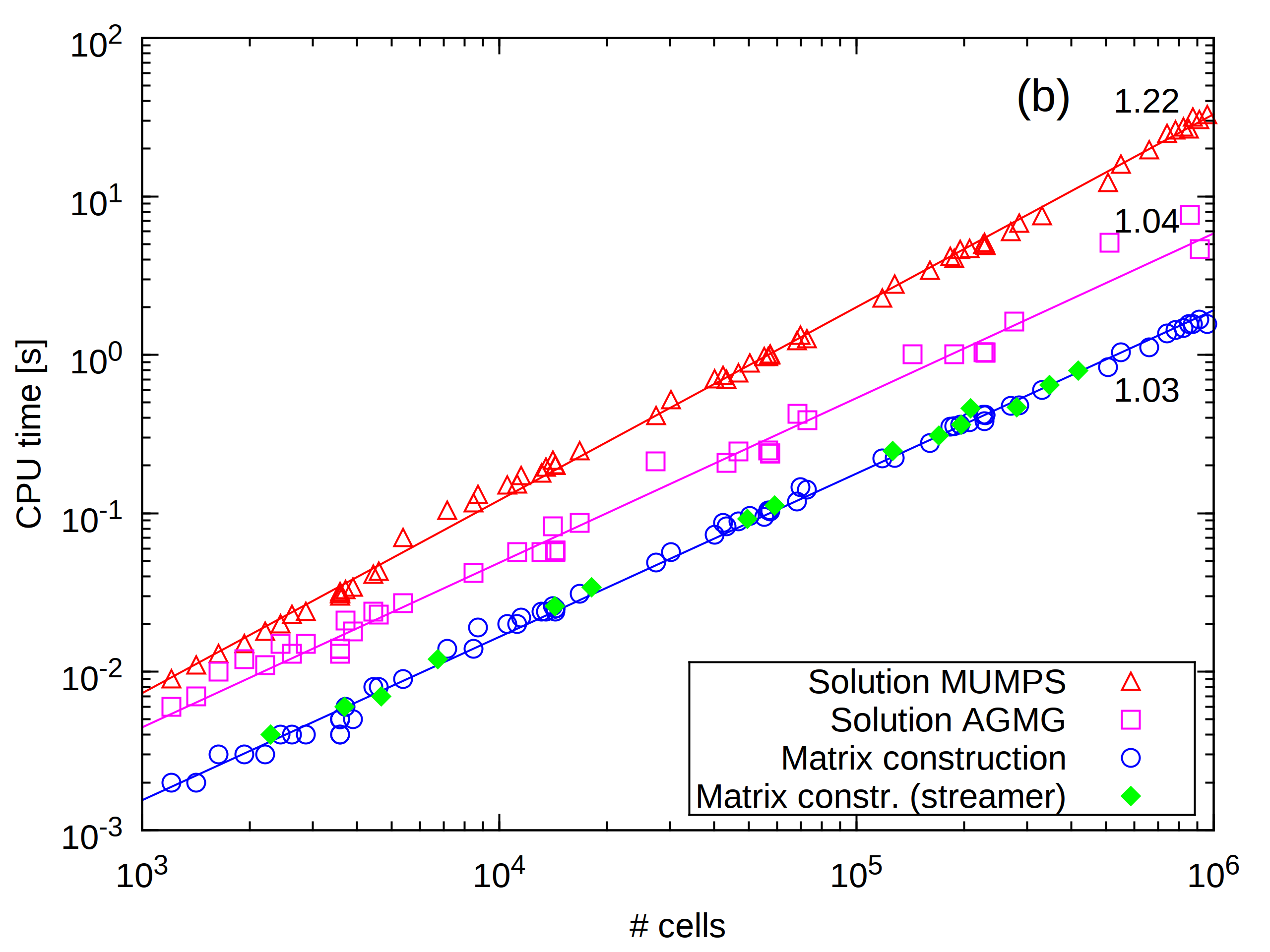}
\caption{Matrix construction and solution: (a) 
discrete Laplacian
$\mathbf{\widetilde{A}}$ on a two-dimensional
multiresolution grid; and (b)
CPU times to build $\mathbf{\widetilde{A}}$
and solve the corresponding linear systems
for several numbers of cells 
(slopes of data fits are indicated).}
\label{fig:matrix} 
\end{figure}  
\begin{figure}[!ht] 
\centering
\includegraphics[width=0.39\hsize]{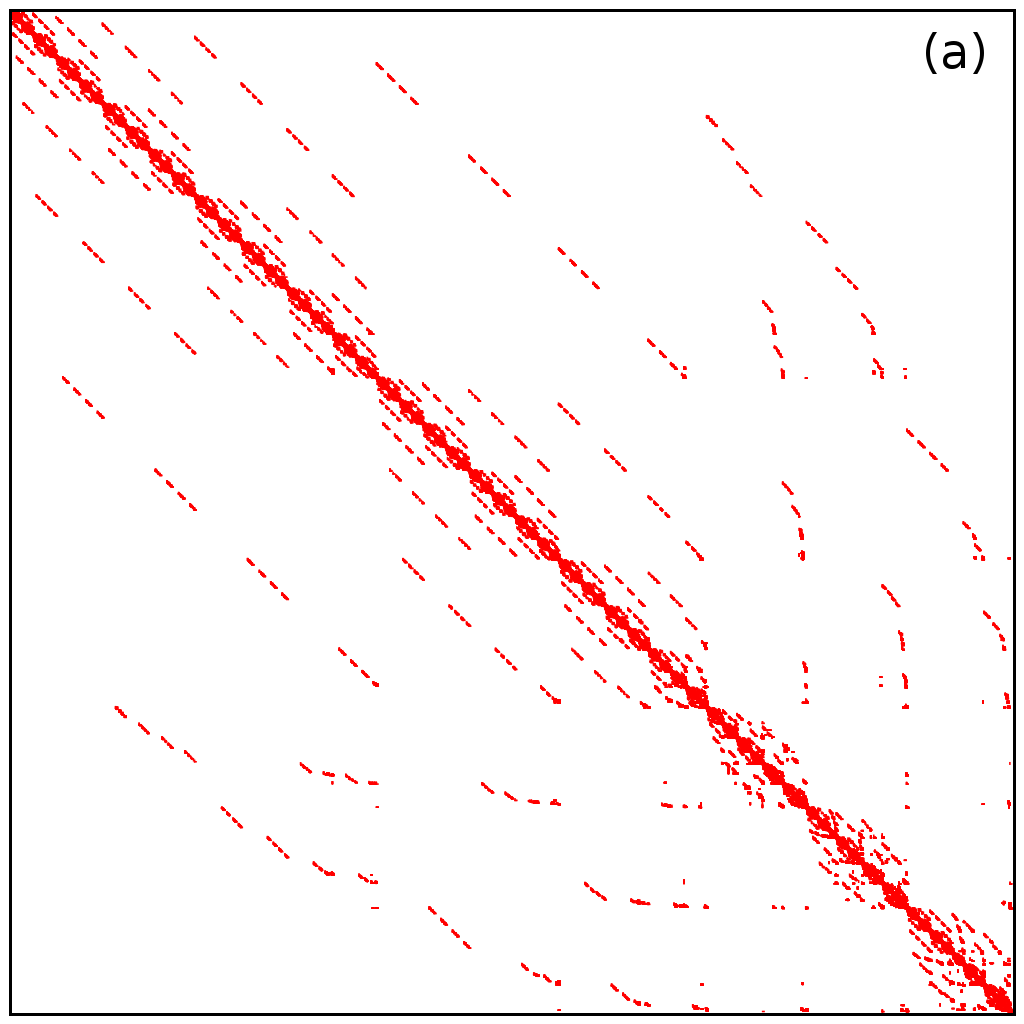}
\includegraphics[width=0.49\hsize]{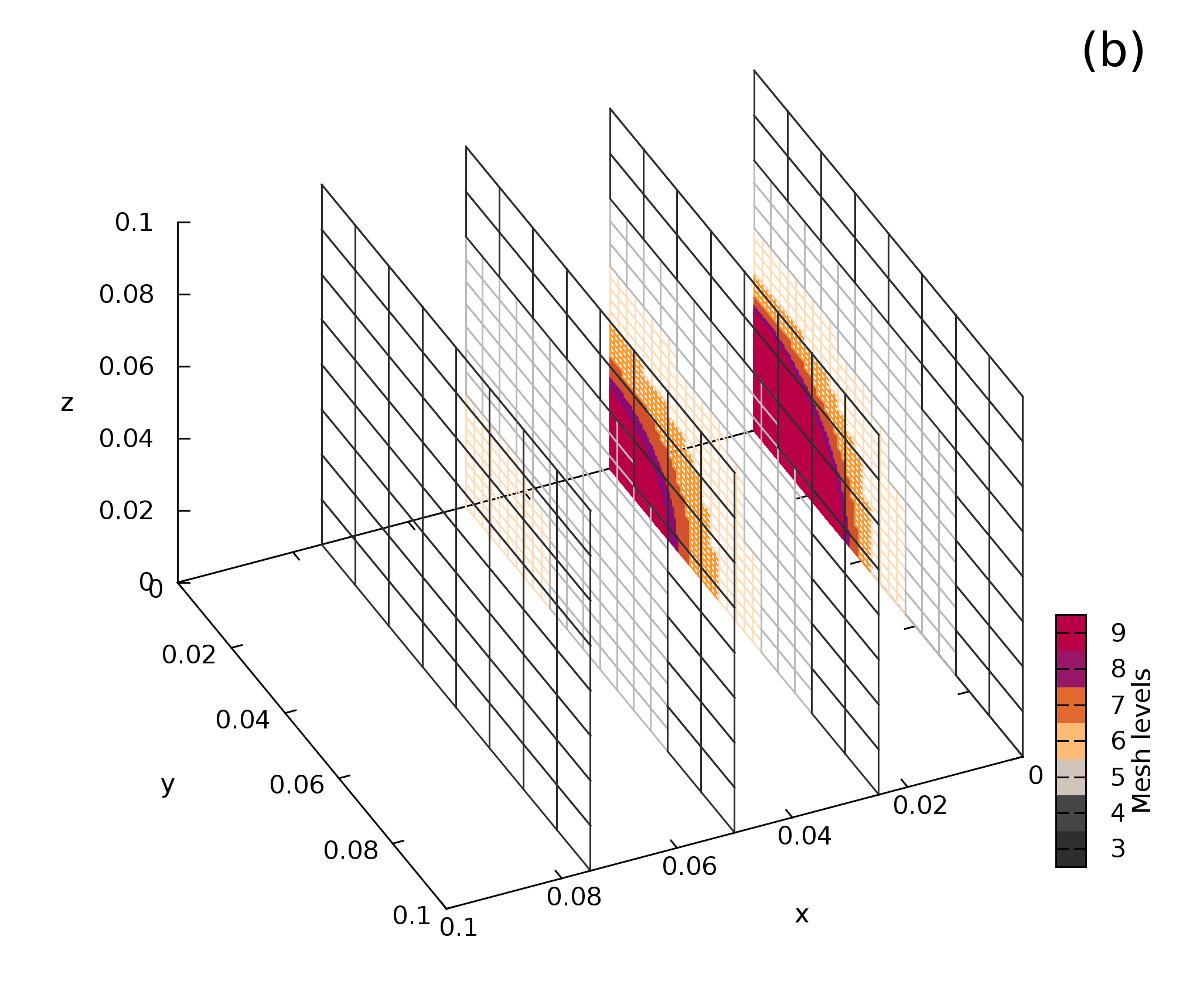}
\caption{Discrete Laplacian
$\mathbf{\widetilde{A}}$ on a 
three-dimensional multiresolution grid (a); and (b) cutting planes 
through the
adapted mesh for $J=9$,
$N_{{\rm R}x}=N_{{\rm R}y}=N_{{\rm R}z}=5$,
and $\eta_{\rm MR}=10^{-4}$.}
\label{fig:matrix3d} 
\end{figure}  

These tests allow us to 
verify that the discrete Laplacian
is consistently constructed following the procedure
established in \S\ref{subsec:const_lap},
and correctly implemented in practice.
The matrix representation $\mathbf{\widetilde{A}}$
is shown in
Figure \ref{fig:matrix}(a)
(recall that only non-zero entries are actually saved in memory).
Notice that matrix $\mathbf{\widetilde{A}}$ cannot be symmetric,
unless no grid adaptation is performed.
However, $\mathbf{\widetilde{A}}$
is in general quasi-symmetric.
For instance,
for this particular problem 
approximately $89\,$\% of symmetry is 
retrieved in terms of non-zero
elements of the matrix.
Figure \ref{fig:matrix}(b) illustrates
the computational complexity of 
the matrix construction, which behaves
like $\Or(N_{\rm L})$.
The different measures were obtained performing
several computations with different finest grid-levels $J\in [9,\ldots,13]$,
and multiresolution parameters $\eta_{\rm MR}\in [10^{-2},\ldots,10^{-9}]$.
We have indicated in Figure \ref{fig:matrix}(b)
the CPU times to solve the resulting linear system with MUMPS, as well as
with an algebraic multigrid solver: AGMG\footnote{
Release 3.1.1.
AGMG home page:
{\tt http://homepages.ulb.ac.be/$\sim$ynotay/AGMG/}} 
\cite{AGMG1,AGMG2,AGMG3}
(tolerances set to $10^{-9}$),
of a computational complexity of $\Or(N_{\rm L})$.
We have also verified that 
building 
the matrix representation behaves the same way
for the more complex streamer configuration
presented in the following.
Finally, Figure \ref{fig:matrix3d}(a) shows operator
$\mathbf{\widetilde{A}}$ on a three-dimensional adapted grid.
As an illustration, the number of non-zero entries is
$29756821$ for $2846787$ cells on the adapted grid 
(a ratio of about $10.5$) with $68\,$\% of symmetry.
The adapted grid for the finest three-dimensional
configuration is shown in Figure \ref{fig:matrix3d}(b),
corresponding to an equivalent uniform grid
of $5^3\times 512^3$ for $J=9$ and
$N_{{\rm R}x}=N_{{\rm R}y}=N_{{\rm R}z}=5$.

\subsection{Performance of linear solvers}
We present a brief study on the performance of several 
software packages currently available in the literature
to solve linear systems of general type: $\mathbf{A}\mathbf{x}=\mathbf{b}$.
Our attention will be focused on iterative solver,
which in general exhibit
relatively modest memory requirements
with respect to direct ones.
In particular algebraic multigrid methods (AMG) 
will be investigated.
These methods
do not require an explicit grid geometry
and work directly on matrix entries;
they are therefore well-suited to our purposes since system
$\mathbf{A}\mathbf{x}=\mathbf{b}$,
stemmed in our case from the discretization of a 
Poisson equation on a multiresolution adapted grid,
has completely lost any reminiscence of its original geometric
layout.

\subsubsection{Test configuration}
Let us consider the propagation of a 
two-dimensional
double-headed streamer at atmospheric pressure.
In this configuration positive and negative streamers emerge
from an initial germ of charged species. 
Drift-diffusion equations (\ref{trasp}) together with Poisson equation
(\ref{poisson})
are solved following the 
time-space
adaptive scheme introduced in \cite{Duarte11_JCP}.
The latter is based on a decoupled numerical solution of (\ref{trasp})
and (\ref{poisson}) in such a way that each problem is solved 
separately by a dedicated solver.
Both numerical approximations are assembled according
to a second order scheme in time.
The latter also considers a time-stepping procedure
with error control such that a prescribed accuracy
$\eta_\Time$ is attained.
Variables are represented 
at cell centers except for the electric field
and the velocities which are staggered,
while the entire problem is solved on an adapted grid
dynamically obtained by multiresolution analysis.
The latter is performed on the 
species densities.
Notice that the right-hand side of the Poisson equation
(\ref{poisson}) is a linear combination of these variables;
hence, the theoretical framework
in \S\ref{subsec:PoissonMR} remains valid.

\begin{figure}[!ht] 
\begin{center}
\includegraphics[width=0.75\hsize]{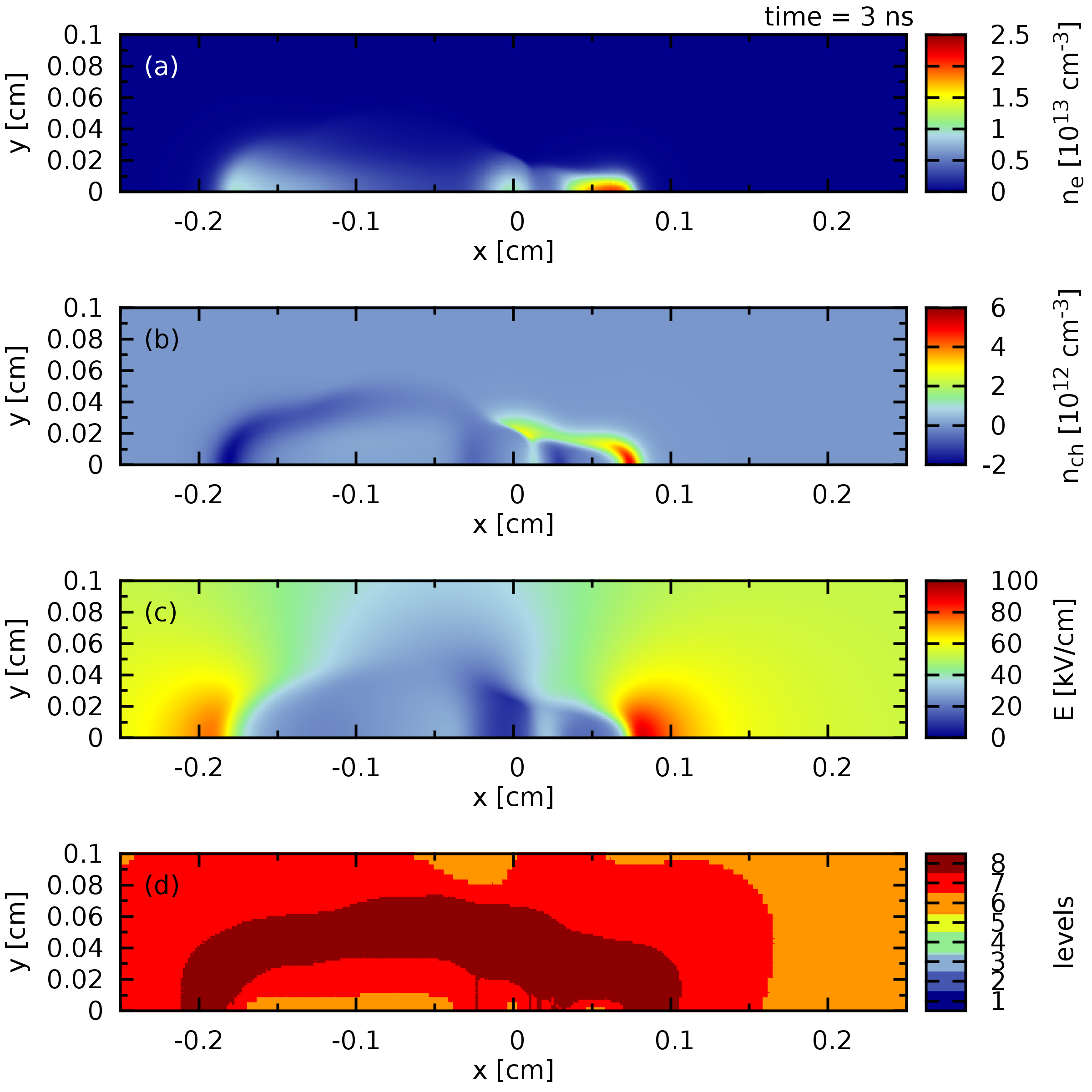}
\end{center}
\caption{Double-headed streamer at $3.0\,$ns: 
         (a) electron density $n_{\rm e}$;
         (b) net charged-species density $n_{\rm ch}$;
         (c) magnitude of the electric field $E=|\vec E|$;
         and
         (d) grid-levels of the adapted mesh. 
Only part of the computational domain is shown.}
\label{DensityPlots} 
\end{figure}

Numerical simulations in the present study were carried out with
a space-time accuracy tolerance of $\eta_{\rm MR} =\eta_\Time = 10^{-4}$
with a space resolution of 3.9$\,\mu$m corresponding to a finest grid level: $J=8$
with $N_{{\rm R}x}=10$ and $N_{{\rm R}y}=3$.
This set of parameters guarantees a sufficiently fine time-space
representation of the physics, and numerical results disclosing practically the same
behavior with higher spatial resolutions and tighter accuracy tolerances.
The computational domain is given by $[-0.5,0.5]\times [0,0.3]\,$cm
in a Cartesian configuration. 
A homogeneous 
electric field $\vec E  = (48.0, 0)\,$kV/cm is
introduced via Dirichlet boundary conditions for the Poisson equation at $x=\pm0.5\,$cm,
whilst Neumann boundary conditions are applied at $y=0.3\,$cm.
A plane of symmetry is imposed  at $y=0$, thus  only one half of the streamer is actually simulated.
The double-headed streamer is initiated by placing a Gaussian plasma cloud so that
the initial conditions for the transport equations (\ref{trasp}) are given by
\begin{equation*}
\np(\x, 0) = \ne(\x, 0) = n_{\max}\exp\left(-|\x|^2/\sigma^2\right)
+ n_{0{\rm p,e}},\quad \nn(\x, 0) = n_{0{\rm n}},
\end{equation*}
with $\sigma=0.02\,$cm, $n_{\max}=10^{13}\,$cm$^{-3}$,
and a small homogeneous pre-ionization background
of $n_{0{\rm n,e}} = 5\times 10^{-5}\,$cm$^{-3}$
and $n_{0{\rm p}} = 10^{-4}\,$cm$^{-3}$.
All tests were conducted starting from the same solution at $3.0\,$ns
when the double-headed streamer is already well developed but no interference 
with the boundaries is evidenced.
The electron density, the net charge species density:
$n_{\rm ch}=\np-\nn-\ne$, the magnitude of the electric field, and
the levels of the adapted grid 
at $3.0\,$ns are shown in
Figure \ref{DensityPlots}.
The total number of cells is of $197784$,
distributed over five grid-levels from a 
resolution of 62.5$\,\mu$m at level $j=4$
up to 3.9$\,\mu$m at $J=8$.
A data compression of about $10\,$\% is thus achieved
with respect to a uniform grid with the finest spatial resolution.
\def\tlvs{\vrule height 1em  width 0pt} 
\begin{table}[!ht]
\centering
\caption{Iterative solvers:
number of iterations (\#iter)
for relative tolerances: $tol=10^{-6},\ldots,10^{-14}$,
CPU computing time, $L^2$-error of $\phi$ 
and $|\vec E|$ with respect to solutions computed with MUMPS,
and memory requirements for each solver.\label{tableSolvers2}
}
\begin{tabular}{clccccccccllllc}
\multicolumn{10}{l}{{\bfseries AGMG}} &
\multicolumn{5}{l}{Memory: 82 MB} \\
\hline
\tlvs& $tol$               &&& \#iter &&& CPU(s) &&& $L^2$-error $\phi$ &&&  $L^2$-error $|\vec E|$ & \\
\hline
\tlvs & $10^{-6}$   &&&       3     &&&    0.50      &&&         1.65$\times10^{-5}$   &&&       1.40$\times10^{-4}$      & \\
\tlvs & $10^{-7}$   &&&       4     &&&    0.55      &&&         1.20$\times10^{-5}$   &&&       2.94$\times10^{-5}$      & \\
\tlvs & $10^{-8}$   &&&       8     &&&    0.78      &&&         1.80$\times10^{-6}$   &&&       4.20$\times10^{-6}$      & \\
\tlvs & $10^{-9}$   &&&      10     &&&    0.89      &&&         1.43$\times10^{-7}$   &&&       4.10$\times10^{-7}$      & \\
\tlvs & $10^{-10}$  &&&      14     &&&    1.10      &&&         2.03$\times10^{-8}$   &&&       4.74$\times10^{-8}$      & \\
\tlvs & $10^{-11}$  &&&      16     &&&    1.21      &&&         2.24$\times10^{-9}$   &&&       7.44$\times10^{-9}$      & \\
\tlvs & $10^{-12}$  &&&      19     &&&    1.38      &&&         2.09$\times10^{-11}$  &&&       9.65$\times10^{-11}$     & \\
\tlvs & $10^{-13}$  &&&      20     &&&    1.43      &&&         1.28$\times10^{-11}$  &&&       3.68$\times10^{-11}$     & \\
\tlvs & $10^{-14}$  &&&      24     &&&    1.64      &&&         1.59$\times10^{-12}$  &&&       3.86$\times10^{-12}$     & \\
\hline
\multicolumn{15}{c}{} \\
\multicolumn{10}{l}{{\bf {\em hypre} BoomerAMG}}  &
\multicolumn{5}{l}{Memory: 100 MB} \\
\hline
\tlvs& $tol$     &&& \#iter &&& CPU(s) &&& $L^2$ error $\phi$ &&&  $L^2$ error $|\vec E|$ & \\
\hline
\tlvs & $10^{-6}$   &&&     3     &&&   1.23     &&&     7.36$\times10^{-4}$                &&&      2.46$\times10^{-3}$     & \\
\tlvs & $10^{-7}$   &&&     6     &&&   1.51     &&&     1.81$\times10^{-5}$                &&&      5.99$\times10^{-5}$     & \\
\tlvs & $10^{-8}$   &&&     8     &&&   1.71     &&&     3.97$\times10^{-6}$                &&&      1.63$\times10^{-5}$     & \\
\tlvs & $10^{-9}$   &&&    11     &&&   1.98     &&&     9.56$\times10^{-8}$                &&&      7.44$\times10^{-7}$     & \\
\tlvs & $10^{-10}$  &&&    14     &&&   2.27     &&&     9.04$\times10^{-9}$                &&&      9.87$\times10^{-8}$     & \\
\tlvs & $10^{-11}$  &&&    17     &&&   2.55     &&&     5.45$\times10^{-10}$               &&&      5.18$\times10^{-9}$     & \\
\tlvs & $10^{-12}$  &&&    20     &&&   2.83     &&&     6.24$\times10^{-11}$               &&&      1.02$\times10^{-9}$     & \\
\tlvs & $10^{-13}$  &&&    24     &&&   3.21     &&&     6.28$\times10^{-12}$               &&&      2.58$\times10^{-11}$    & \\
\tlvs & $10^{-14}$  &&&    27     &&&   3.52     &&&     4.73$\times10^{-13}$               &&&      3.84$\times10^{-12}$    & \\
\hline
\multicolumn{15}{c}{} \\
\multicolumn{10}{l}{{\bf {\em hypre} BoomerAMG + GMRES} }  &
\multicolumn{5}{l}{Memory: 146 MB} \\
\hline
\tlvs& $tol$     &&& \#iter &&& CPU(s) &&& $L^2$ error $\phi$ &&&  $L^2$ error $|\vec E|$ & \\
\hline
\tlvs & $10^{-6}$   &&&      2     &&&   1.24     &&&   9.09$\times10^{-4}$      &&&    2.36$\times10^{-3}$               & \\
\tlvs & $10^{-7}$   &&&      5     &&&   1.57     &&&   2.65$\times10^{-5}$      &&&    1.28$\times10^{-4}$               & \\
\tlvs & $10^{-8}$   &&&      8     &&&   1.90     &&&   1.15$\times10^{-6}$      &&&    1.49$\times10^{-5}$               & \\
\tlvs & $10^{-9}$   &&&     10     &&&   2.13     &&&   6.19$\times10^{-8}$      &&&    8.56$\times10^{-7}$               & \\
\tlvs & $10^{-10}$  &&&     12     &&&   2.34     &&&   4.32$\times10^{-9}$      &&&    6.09$\times10^{-8}$               & \\
\tlvs & $10^{-11}$  &&&     14     &&&   2.58     &&&   5.77$\times10^{-10}$     &&&    3.33$\times10^{-9}$               & \\
\tlvs & $10^{-12}$  &&&     15     &&&   2.69     &&&   3.58$\times10^{-10}$     &&&    8.20$\times10^{-10}$              & \\
\tlvs & $10^{-13}$  &&&     17     &&&   2.93     &&&   3.57$\times10^{-11}$     &&&    7.10$\times10^{-11}$              & \\
\tlvs & $10^{-14}$  &&&     19     &&&   3.15     &&&   3.13$\times10^{-12}$     &&&    7.73$\times10^{-12}$              & \\
\hline
\end{tabular}
\end{table}

\subsubsection{Analysis of results}
We have considered some iterative
solvers readily available in various software packages.
Most of present day linear solvers are developed with a special attention 
on enhanced parallel capabilities. 
Nevertheless, thanks to significant data compression
achieved by multiresolution adaptation,
the linear systems under consideration have typically 
about $10^5$ unknowns with
approximately $10^6$ non-zero elements in the system matrix.
Therefore, to simplify our study 
we have focused our attention on sequential
performance of these solvers.
We have performed the numerical experiments on a two-processor computer.
Each processor is an Intel Xeon CPU E5410 @ 2.33GHz with a 
total available computer memory of 24 GB. 
The computer runs on a 64-bit version of Fedora 18 GNU/Linux system.
All codes with the various linear solvers were compiled using compilers from GCC (version 4.7.2).
Memory requirements of each solver were obtained by tracing
the memory profiles of running programs with {\tt top} command, executed in batch mode 
with a delay-time 
interval set to $0.01\,$s. 
In order to discriminate memory requirements
for the linear solvers from the overall program memory usage, 
a reference program was executed in which calls to the  
solver were replaced by FORTRAN (GNU extension) {\tt SLEEP} command.

The total number of unknowns for the Poisson equations 
considered in this problem is
given by the number of cells in the adapted grid,
$197784$ in this case, while the discrete
Laplacian has $1078534$ non-zero entries
(a ratio of about $5.5$).
In what follows we consider as reference solution 
the solution to the Poisson equation (\ref{poisson}):
$\phi$,
computed with MUMPS.
For this problem, MUMPS requires $193\,$MB of memory space
for a computation that takes approximately $3.96\,$s.
As before we also analyze the approximation to the electric field:
$\vec E =- \Dx \phi$.
Data for three iterative solvers are presented in 
Table \ref{tableSolvers2} for two algebraic multigrid solvers:
AGMG and BoomerAMG, 
and for GMRES \cite{GMRES} preconditioned with 
BoomerAMG (also contained in \hypre).
In all cases a fine-tuning of computing parameters
have been previously carried out so that
Table \ref{tableSolvers2} includes the best performances
obtained with each of these solvers for this particular
problem.
A key parameter for iterative solvers is given by the
relative and absolute tolerances that in particular
serve as stopping criteria
to the iterative procedures.
In this study we have set both tolerances equal to an
accuracy tolerance, denoted as $tol$.
The initial guess corresponds to the solution computed
during the previous time-step.
For tolerances higher or equal to $10^{-5}$
convergence is attained right-away
with the initial guess for all three solvers.
In all cases better performances are obtained 
with these iterative solvers with respect to MUMPS
even with very fine accuracy tolerances $tol$.
Even though GMRES converges in a less number of iterations
for different values of $tol$ 
with respect to the
algebraic multigrid solvers,
it does not yield faster computations taking
into account that for this problem 
preconditioning is the most expensive 
part.
Therefore,
BoomerAMG and GMRES/BoomerAMG
require similar computing times.

With respect to the reference solution, all
these iterative solvers scale well in terms of the
accuracy of the approximations, set by the
tolerance parameter $tol$.
Notice that numerical errors related to iterative solvers
must be taken into account 
to track the overall numerical accuracy of the simulation.
In particular these
numerical errors must be smaller than the multiresolution
ones so that (\ref{eq:res_coro}) and (\ref{eq:hat_tilde})
remain valid.
The latter could be enforced by setting
in general: $tol < \eta_{\rm MR}$, while in this particular case
a safer choice might be given by $tol \leq 10^{-3}\times \eta_{\rm MR}$
according to the values contained in Table \ref{tableSolvers2}.
Among the solvers tested in this study, AGMG
revealed itself as
the most performing package both in terms of CPU time and memory requirements
to solve this particular problem.
However, the overall performance 
of these solvers is clearly problem-dependent.
In this regard the \hypre library
provides a user-friendly and unified interface to various solution schemes,
very appropriate to handle different types of problems.

\subsection{Application to the study of two interacting positive streamers}
\label{CaseStudy}
While previous illustrations served to validate
the numerical strategy, 
we consider now an interesting plasma physics application
with more complex dynamics.
We study 
the interaction 
of two positive streamers initiated to develop side by side.
Because the heads of both streamers carry space charge of the same polarity,
their mutual interaction should essentially be an electrostatic repulsion.
However, it was found that streamers in such a configuration may attract 
each other and eventually merge \cite{Briels2006,Cummer2006,Nijdam2009}.
This attraction is mainly the result of the enhancement of photoionization
source in the space between the streamer heads \cite{Luque2008,Bonaventura2012}.
In particular, 
based on an extensive parametric numerical study, we have
shown in  \cite{Bonaventura2012}
that for initial separations of two streamers smaller or comparable
to the absorption length of photoionization, merging will start 
when the ratio of the streamer characteristic width and 
their mutual separation attains a certain value.
We describe here some
numerical aspects omitted in our previous study.
\begin{figure}[!ht]
\centering
\includegraphics[width=0.6\hsize]{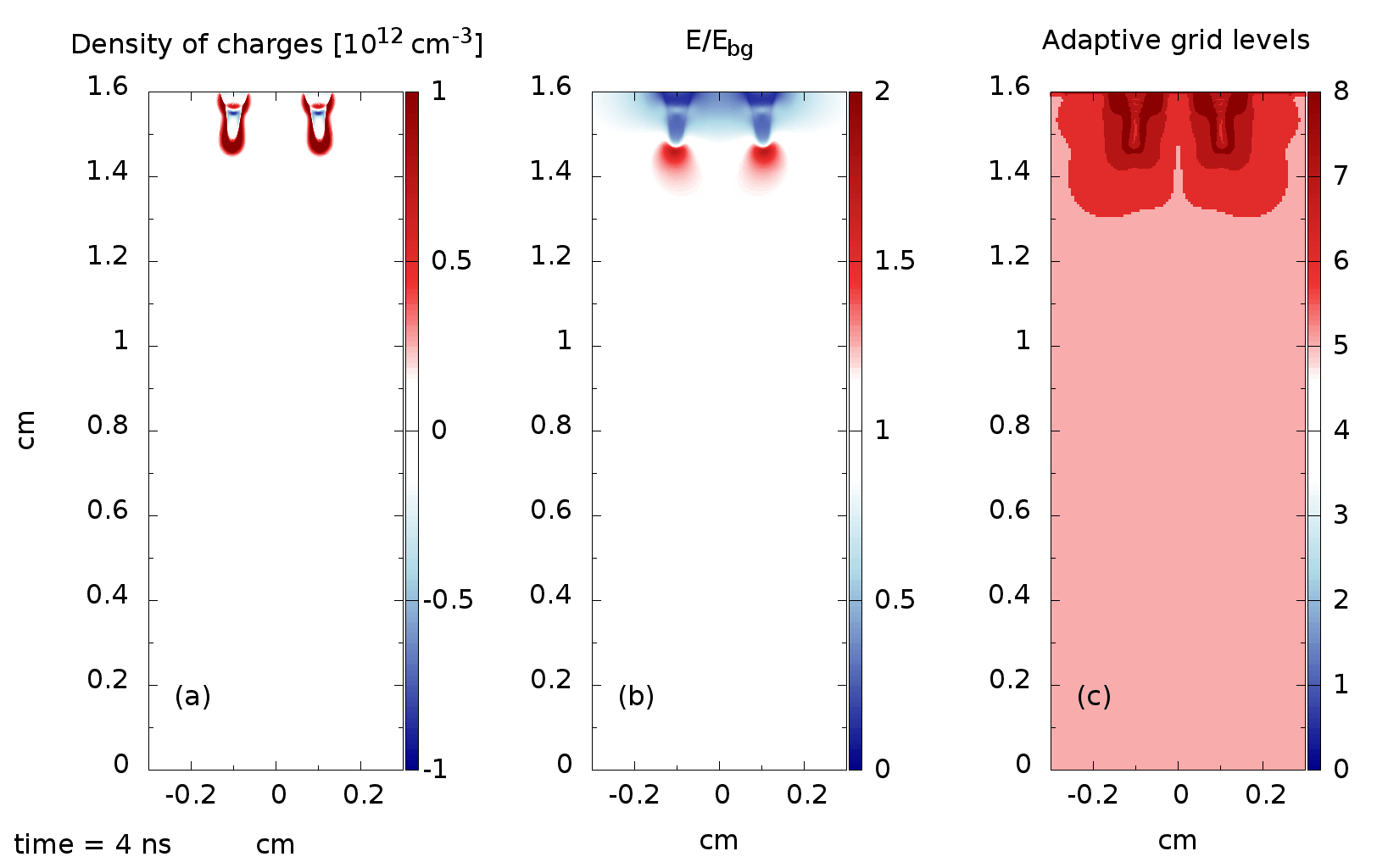}
\rule{0.5\textwidth}{0.2mm}
\includegraphics[width=0.6\hsize]{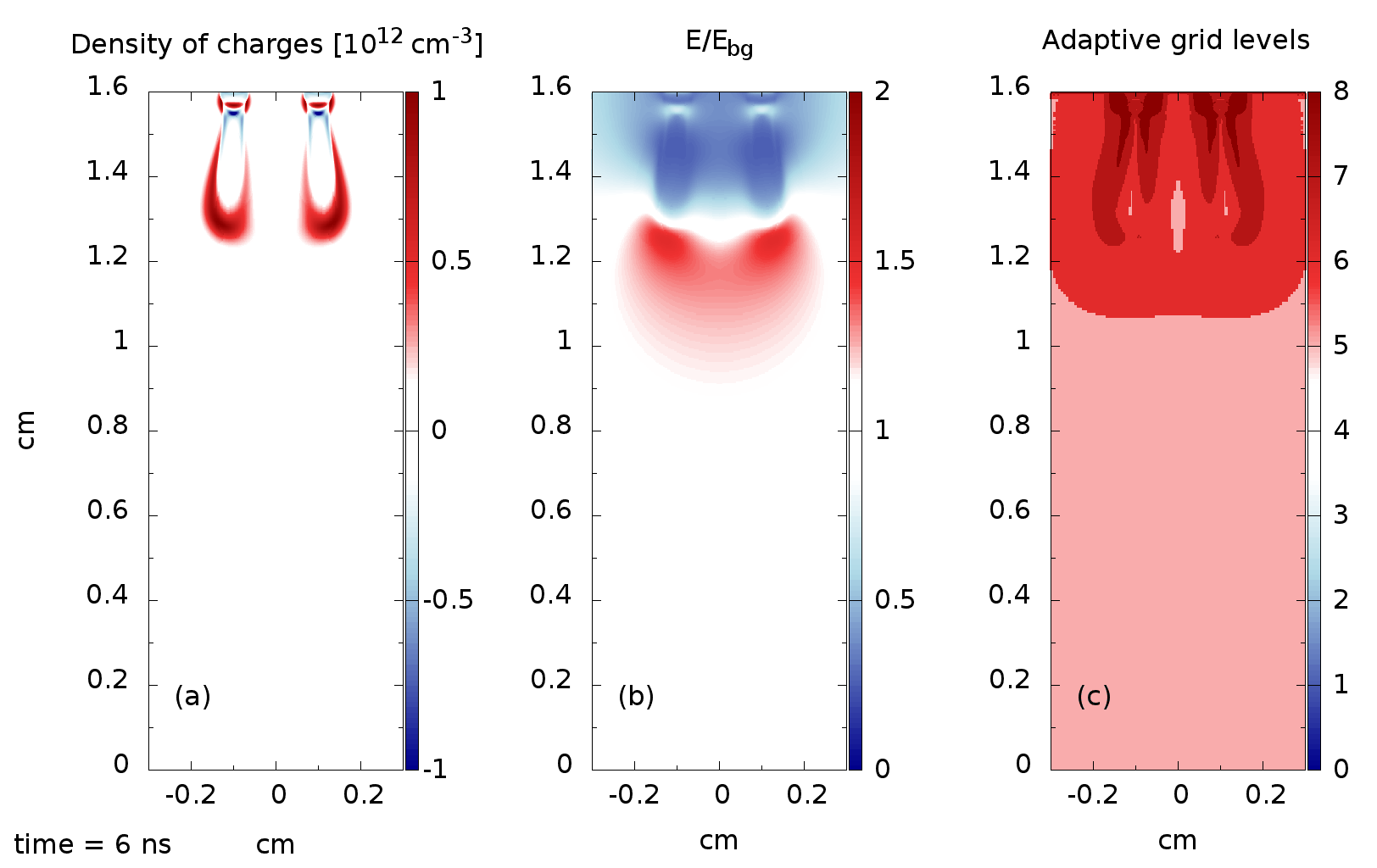}
\rule{0.5\textwidth}{0.2mm}
\includegraphics[width=0.6\hsize]{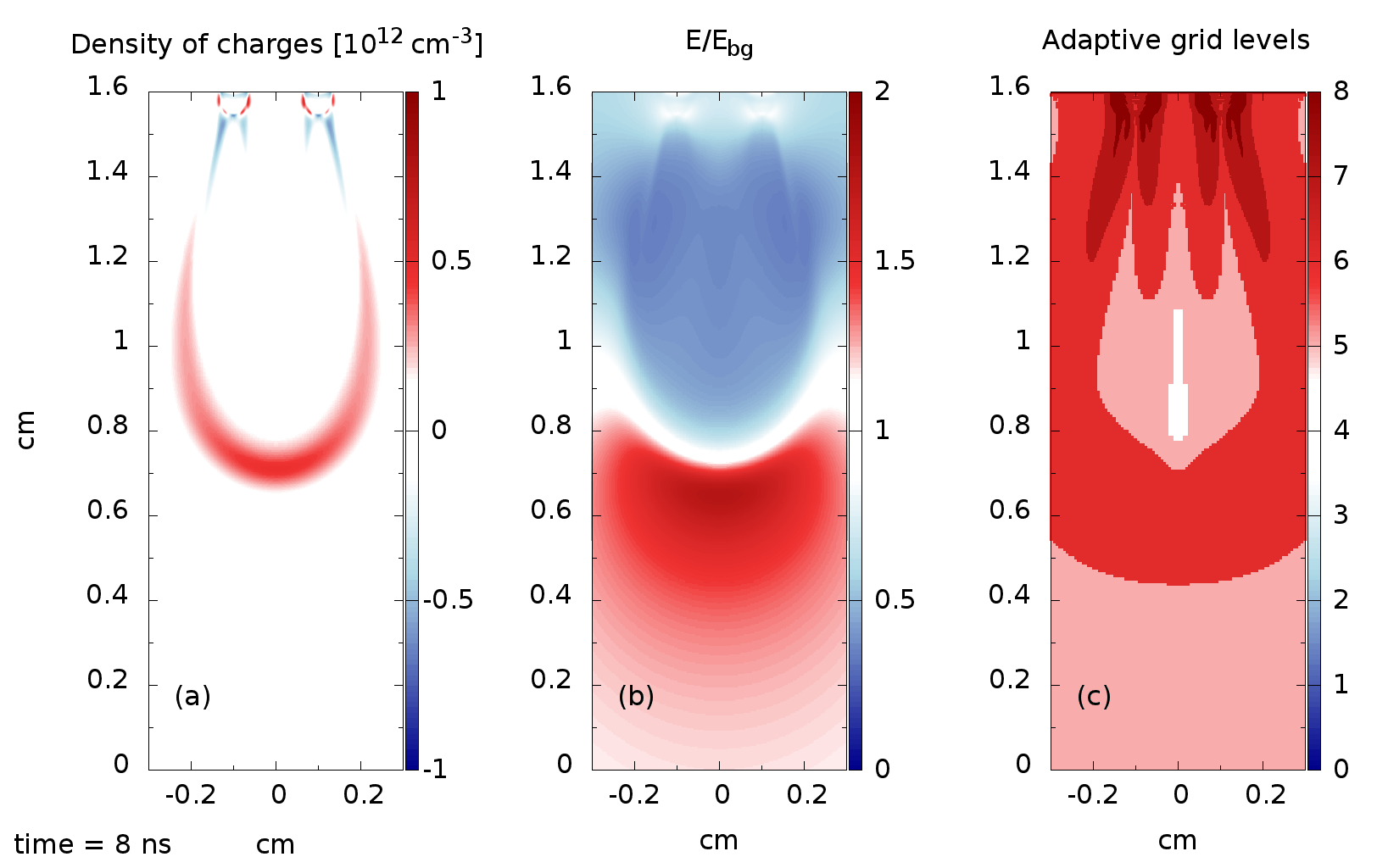}
\caption{Time evolution of the net charge density (a),
magnitude of the electric field (b), and 
dynamic grid adaptation (c) for two interacting positive streamers
at ground pressure  with an applied electric field
of $\vec E_{\rm bg}=(0,-48)\,$kV/cm at time instances: $4.0$, $6.0$ and $8.0\,$ns.
Only part of the computational domain is shown.
}
\label{rho-phases}
\end{figure}

Let us consider two positive streamers
modeled by (\ref{trasp})--(\ref{poisson})
propagating in a homogeneous electric field.
As before the system of equations is solved with 
the time-space adaptive scheme introduced in \cite{Duarte11_JCP}
with the Poisson equation discretized on the adapted grid following
the numerical technique established in \S\ref{subsec:const_lap}.
The resulting linear systems are solved with MUMPS.
The computational domain is given by $[0,3.0]\times[-1.6,1.6]\,$cm
in a Cartesian configuration.
A space-time accuracy tolerance of $\eta_{\rm MR} =\eta_\Time = 10^{-4}$
was chosen 
with a space resolution of 3.9$\,\mu$m corresponding to a finest grid-level of $J=8$
with $N_{{\rm R}x}=30$ and $N_{{\rm R}y}=32$. 
The finest spatial resolution is equivalent to that of a uniform grid
with $8192\times7680$ cells.	
A homogeneous electric field  of $\vec E_{\rm bg} = (0,-48)\,$kV/cm is
introduced via Dirichlet boundary conditions for the Poisson equation (\ref{poisson})
at $y=\pm1.6\,$cm, whilst Neumann boundary conditions are applied at $x=3.0\,$cm.
A plane of symmetry is imposed  at $x=0$. 
The positive streamer is initiated by placing a Gaussian seed
with a maximum  of $10^{13}\,$cm$^{-3}$ and a characteristic
width of $0.02 \,$cm, centered at $0.1\,$cm from the symmetry axis. 
The time evolution of the net charge density, the magnitude of the electric field,
and the dynamic grid adaptation at time instances: $4.0$, $6.0$, and $8.0\,$ns
are shown in Figure~\ref{rho-phases}.
Population of different grid-levels at sample times is detailed in
Table~\ref{TableResol} together with the corresponding data compression (DC),
defined as the percentage of active cells with respect to the equivalent 
number of cells using the finest discretization,
given in this case by 62914560.
We recall that no grid overlapping is considered in this implementation,
that is, both the time-dependent PDEs as well as the Poisson equations
are solved on the adapted grid consisting of cells at different grid-levels
as shown in Table~\ref{TableResol}.
The coarsest resolution allowed in this simulation (at grid-level $j=1$)
corresponds to a spatial resolution of $0.05\,$cm (note that
this level was not populated during the simulation,
therefore it is not listed in Table~\ref{TableResol}).
\begin{table}[!ht]
\caption{Data compression (DC) and number of cells at different grid-levels at sample time instances.}
\label{TableResol}
\begin{center}
\begin{tabular}{ccrrrrrrr} 
\toprule
\multicolumn{1}{c}{time(ns)}  & \multicolumn{1}{c}{DC(\%)} &  \multicolumn{7}{c}{Number of cells at grid-levels} \cr
& 
      &  \multicolumn{1}{c}{2}   
      &  \multicolumn{1}{c}{3}
      &  \multicolumn{1}{c}{4}  
      &  \multicolumn{1}{c}{5}  
      &  \multicolumn{1}{c}{6}
      &  \multicolumn{1}{c}{7}   
      &  \multicolumn{1}{c}{8}  \cr
\midrule
\phantom{0}0  &    0.031  & 15240 & 162  &  272 &   527 &  1341 & 7728  &169904\cr
\phantom{0}1  &    1.219  & 780    &4671&61359&605154&22297&23160&49424 \cr
\phantom{0}2  &    1.264  & 620  & 4440 &57625 &633483 & 23447 &27668  & 47744\cr
\phantom{0}3  &    1.291  & 566  &4346&56348&642443&26479&33016&48752\cr 
\phantom{0}4  &    1.325  & 528  & 4528 &54141 &648333 & 34030 &41660  & 50672\cr 
\phantom{0}5  &    1.360  & 512 & 4712 & 51905 & 652364 & 42504 & 56536 & 47136 \cr
\phantom{0}6  &    1.407  & 544  & 4768 &49464 &653958 & 61267 &64120  & 51056\cr
\phantom{0}7  &    1.449  & 458 & 5288 & 47652 & 650502 & 93008 & 63196 & 51728 \cr
\phantom{0}8  &    1.495  & 338  & 5791 &46530 &644986 &131087 &66656  & 45328\cr
\phantom{0}9  &    1.567  & 219 & 5388 & 48316 & 637538 & 188227 & 70380 & 35744 \cr
10 &    1.702  & 18 & 5296 & 48622 & 624126 & 294256 & 73264 & 25280 \cr 
\bottomrule
\end{tabular}
\end{center}
\end{table}

From Figure~\ref{rho-phases} and Table~\ref{TableResol}
we observe that the finest level is first populated at the vicinity 
of the initial Gaussian seed and follows the propagation of onsetting streamers
(see Figure~\ref{rho-phases} corresponding to time $4\,$ns). 
At the next instance shown (Figure~\ref{rho-phases} at $6\,$ns),
the propagating  front is fully  described in a region 
contained within levels 6 and 7. This is
because both streamer heads had expanded and the finest scale
is thus no longer necessary.
Once the streamer heads have merged (Figure~\ref{rho-phases} at $8\,$ns)
and therefore only one head is propagating, 
only level 6 is required. 
In particular
behind the
head, {\it i.e.}, inside the plasma channel where neither sharp gradients
nor strong discharge activity are present, 
the grid is coarsened down to level 4.
The finest resolution is attained and kept 
throughout the simulation
close to the initial Gaussian seeds where we can observe
persistence of highly localized space charge  
as well as strong spatial variation of the electric field.
Despite a decreasing population of the finest level after $7\,$ns
(see Table~\ref{TableResol})
overall data compression is slowly increasing
because discharge activity
is gradually filling larger regions of the computational domain.

\section{Concluding remarks}\label{sec:conclusion}
The multiresolution finite volume scheme \cite{Cohen03}
has been extended to include the numerical solution of Poisson equations
on the corresponding adapted grids.
A numerical procedure
has been developed
to represent the discrete Laplace operator
on the adapted grid
by reconstructing locally uniform-grid regions
at inter-grid interfaces by means of ghost cells
and inter-level multiresolution operations.
This approach constitutes a new alternative
to the standard level-wise numerical solution of
elliptic equations considered in most of the 
adaptive mesh refinement techniques 
for time-dependent problems
in the literature.
The numerical solution of 
the discrete Poisson equation amounts to considering a linear
system completely independent of the grid 
generation or any other grid-related
data structure or geometric consideration.
The multiresolution framework
guarantees numerical approximations within an accuracy
tolerance as well as consistency and conservation
properties throughout the set of grids.
Here we have focused our attention on Poisson
equations, however the present technique 
remains valid for more general elliptic PDEs
like Poisson equations with
time- and/or space-varying
coefficients.

The validity of 
the numerical strategy
has been assessed
in the context of the numerical simulation of 
streamer discharges.
This application involves an intensive
use of Poisson solvers and accurate
solutions of Poisson equations are essential to
the correct reproduction of physics.
First,
we have carefully evaluated the 
numerical errors introduced by data compression
for a
simpler configuration with analytical solution.
A much more complex and complete model
was then considered to simulate the propagation of a double-headed streamer
discharge in air at atmospheric pressure.
We have thus conducted a study on the performance
and capabilities of various 
linear solvers for this problem.
The latter allowed us
to further validate the current implementation and 
serves as a guide for other applications.
In particular we have evaluated the potentialities of
algebraic multigrid solvers, well-suited for this 
kind of implementation with no geometric counterpart.
The robustness of the numerical strategy has been further
assessed for the simulation of 
interacting positive streamers,
an interesting application in  
plasma physics.

Further developments include optimizing  
the numerical construction of the discrete Laplace
operators by conceiving, for instance, better
data structures
or 
by updating only the matrix entries modified by
grid adaptation.
Taking into account that in this implementation
solving the linear systems
becomes a separate aspect from the
multiresolution analysis itself,
parallel computing capabilities
may be directly inherited from the software packages
available in the literature.
However,
an intelligent conjunction with multiresolution
parallelism
must be sought to achieve overall satisfactory results.
These issues constitute particular topics of our
current research.

\section*{Acknowledgments}
This research was 
supported by a fundamental project grant from ANR (French National Research Agency - ANR Blancs):
{\it S\'echelles} (project leader S.~Descombes - 2009-2013)
and 
by a DIGITEO RTRA project: \emph{MUSE} (project leader M.~Massot - 2010-2014).
M.~D. acknowledges support of 
Laboratoire EM2C for a visiting stay in France.
Z.~B. acknowledges support from project CZ.1.05/2.1.00/03.0086 
funded by the European Regional Development Fund and support of Ecole Centrale Paris.

\appendix

\section{Multiresolution error estimate}\label{AppMR}
Defining the pairs $(\phi_{j,k},\psi_{j,k})$ and $(\widetilde{\phi}_{j,k},\widetilde{\psi}_{j,k})$
as, respectively, the {\it primal} and the {\it dual}
scaling function and wavelet,
the following representations of $f_{j+1}$
are perfectly equivalent:
\begin{equation*}\label{eq3:f_j_interlevel}
f_{j+1} :=  \ds \sum_{|\mu|=j+1 } f_\mu \phi_\mu = 
\ds \sum_{|\gamma|=j } f_\gamma \phi_\gamma +
\ds \sum_{|\gamma|=j} d_\gamma \psi_\gamma,
\end{equation*}
with 
$f_\mu:=\langle f,\widetilde{\phi}_\mu\rangle$,
$f_\gamma:=\langle f,\widetilde{\phi}_\gamma\rangle$, and
$d_\gamma:=\langle f,\widetilde{\psi}_\gamma\rangle$.
Iterating on $j$, we have the following wavelet decomposition,
\begin{equation*}
 f = \ds \sum_{j = -1}^{\infi} \sum_{|\lambda|=j} \langle f,\widetilde{\psi}_{\lambda}\rangle \psi_{\lambda},
\end{equation*}
where  $\psi_{-1,k}:=\phi_{0,k}$ and  $\widetilde{\psi}_{-1,k}:=\widetilde{\phi}_{0,k}$.
We can thus construct
the array $\mathbf{\Psi}_{J,\lambda}$,
$|\lambda|\leq J$, that corresponds to 
the primal wavelets $\psi_{\lambda}$
cell-averaged at level $J$, {\it i.e.},
$\mathbf{\Psi}_{J,\lambda} := (\langle \psi_{\lambda},\widetilde{\phi}_{\gamma}\rangle)_{\gamma \in S_J}$.
For compactly supported wavelets,
there is a constant $C>0$ such that 
\begin{equation*}
 \|\mathbf{\Psi}_{J,\lambda}\|_2 \leq C \|\psi_{\lambda} \|_{L^2} \leq C 2^{-d|\lambda|/2},
\end{equation*}
and for the multiresolution approximation $\adap_{\Lambda}\mathbf{f}_J$, 
we have that
\begin{equation*}\label{eq3:error_adap}
\|\mathbf{f}_{J} - \adap_{\Lambda}\mathbf{f}_J \|^2_2 = 
\left\| d_{\lambda} \mathbf{\Psi}_{J,\lambda}\vert_{\lambda \notin \Lambda} \right\|^2_2
\leq C \ds \sum_{\lambda \notin \Lambda} \|d_{\lambda}\|_{L^2}^2 2^{-d|\lambda|} =
C \ds \sum_{\|d_{\lambda}\|_{L^2} \leq \epsilon_{|\lambda|}} 
\|d_{\lambda}\|_{L^2}^2 2^{-d|\lambda|},
\end{equation*}
because only some of the components of $(\mathbf{f}_{J} - \adap_{\Lambda}\mathbf{f}_J)$
are non-zero, namely those corresponding to discarded details;
therefore, the approximation error is bounded by their sum.
Considering a level-wise threshold parameter:
$\epsilon_j := 2^{dj/2} \epsilon$,
the next bound follows (where $\#(\cdot)$ returns the cardinality of a set):
\begin{equation*}\label{eq3:error_adap_2}
\|\mathbf{f}_{J} - \adap_{\Lambda}\mathbf{f}_J \|^2_2 \leq
C \#(\nabla^J)\epsilon^2 = C \# (S_J)\epsilon^2 \leq C 2^{dJ}\epsilon^2,
\end{equation*}
with the cautious assumption that 
$\|d_{\lambda}\|_{L^2}=\epsilon_{|\lambda|}$
for all $d_{\lambda}$ such that   
$\lambda \notin \Lambda$
(even though they might be much smaller than $\epsilon_{|\lambda|}$)
as well as for the remaining components of 
$(\mathbf{f}_{J} - \adap_{\Lambda}\mathbf{f}_J)$
(even though they are zero).
Choosing $\epsilon := 2^{-dJ/2}\eta_{\rm MR}$
then yields (\ref{eq3:adap_error_eps})
with the level-dependent threshold values (\ref{eq3:epsilon_j}).
Bound (\ref{eq3:adap_error_eps}) is similarly shown in \cite{Cohen03} 
for both a uniform and $\ell^1$ norms.

\section{Pseudo-code of the algorithm}\label{AppScheme}
We consider 
a multiresolution
adapted grid
given by the set of leaves:
$\Theta_\leaf= (\Omega_{\lambda})_{h(\lambda)\in {\rm I}_\leaf }$,
a one-dimensional array
of size $N_{\rm L}$.
The algorithm to construct the discrete
Laplacian: 
$\mathbf{\widetilde{A}} \in \M$,
can be schematically described as follows
in a Cartesian finite volume 
framework where interfaces are given by
$\Gamma^{d'}_{\gamma,\mu}$, $d'=1,\ldots,d$.
This scheme supports 
polynomial interpolations (\ref{eq3:linear_prediction})
and
finite volume space discretizations
(\ref{eq:flux})
of arbitrary order.
\begin{algorithmic}
\STATE {\bf Initialization:} $\mathbf{\widetilde{A}}=0$. 
\FOR{$i = 1 \to N_{\rm L} $}
\STATE Current leaf: $\Omega_\gamma$ s.t. $\gamma = h^{-1}(i)$.
\FOR{$d' = 1 \to d$}
\STATE Current neighbor: 
$\Omega_{\mu}$ s.t. 
$\Gamma^{d'}_{\gamma,\mu}=\overline{\Omega_{\gamma}} \cap \overline{\Omega_{\mu}}$.
\IF[$\Omega_\mu$ is a leaf, {\it i.e.}, (\ref{eq:leaf_a}).]{$\mu \in D(h)$ }
\STATE $i'=h(\mu)$.
\FOR{$\lambda \in R^+_F(\gamma)$}
\IF[$\Omega_\lambda$ is a leaf, {\it i.e.}, (\ref{eq:leaf_a}).]{$\lambda \in D(h)$ }
\STATE $l=h(\lambda)$.
\STATE $\widetilde{a}_{i,l} = \widetilde{a}_{i,l} + \alpha_{\gamma,\lambda}$. 
\STATE $\widetilde{a}_{i',l} = \widetilde{a}_{i',l} - \alpha_{\gamma,\lambda}$.
\ELSIF[$\Omega_\lambda$ is a phantom, {\it i.e.}, (\ref{eq:phantom_a}).]{$\lambda \in \bigcup_{r=1}^{N_{\rm R}}  \phan(\Lambda_r)$ }
\FOR{$\widehat{\lambda}$ s.t. $\widehat{\lambda} \in R_I(\lambda)$}
\IF[$\Omega_{\widehat{\lambda}}$ is a leaf, {\it i.e.}, (\ref{eq:leaf_a}).]{$\widehat{\lambda} \in D(h)$ }
\STATE $l=h(\widehat{\lambda})$.
\STATE $\widetilde{a}_{i,l} = \widetilde{a}_{i,l} + \beta_{\lambda,\widehat{\lambda}}\alpha_{\gamma,\lambda}$.
\STATE $\widetilde{a}_{i',l} = \widetilde{a}_{i',l} - \beta_{\lambda,\widehat{\lambda}}\alpha_{\gamma,\lambda}$.
\ELSE[$\Omega_{\widehat{\lambda}}$ is within the tree, {\it i.e.}, (\ref{eq:tree_a}).]
\FOR{$\widehat{\lambda}'$ s.t. $\Omega_{\widehat{\lambda}'} \subset \Omega_{\widehat{\lambda}}$}
\IF[$\Omega_{\widehat{\lambda}'}$ is a leaf, {\it i.e.}, (\ref{eq:leaf_a}).]{$\widehat{\lambda}' \in D(h)$ }
\STATE $l=h(\widehat{\lambda}')$.
\STATE $\widetilde{a}_{i,l} = \widetilde{a}_{i,l} + 
|\Omega_{\widehat{\lambda}}|^{-1}|\Omega_{\widehat{\lambda}'}| \beta_{\lambda,\widehat{\lambda}}\alpha_{\gamma,\lambda}$.
\STATE $\widetilde{a}_{i',l} = \widetilde{a}_{i',l} - 
|\Omega_{\widehat{\lambda}}|^{-1}|\Omega_{\widehat{\lambda}'}| \beta_{\lambda,\widehat{\lambda}}\alpha_{\gamma,\lambda}$.
\ELSE
\FOR{$\widehat{\lambda}''$ s.t. $\Omega_{\widehat{\lambda}''} \subset \Omega_{\widehat{\lambda}'}$}
\IF[$\Omega_{\widehat{\lambda}''}$ is a leaf, {\it i.e.}, (\ref{eq:leaf_a}).]{$\widehat{\lambda}'' \in D(h)$ }
\STATE $l=h(\widehat{\lambda}'')$.
\STATE $\widetilde{a}_{i,l} = \widetilde{a}_{i,l} + 
|\Omega_{\widehat{\lambda}}|^{-1}|\Omega_{\widehat{\lambda}''}|\beta_{\lambda,\widehat{\lambda}}\alpha_{\gamma,\lambda}$.
\STATE $\widetilde{a}_{i',l} = \widetilde{a}_{i',l} - 
|\Omega_{\widehat{\lambda}}|^{-1}|\Omega_{\widehat{\lambda}''}|\beta_{\lambda,\widehat{\lambda}}\alpha_{\gamma,\lambda}$.
\ELSE
\STATE Continue up to leaves.
\ENDIF
\ENDFOR
\ENDIF
\ENDFOR
\ENDIF
\ENDFOR
\ELSE[$\Omega_{\lambda}$ is within the tree, {\it i.e.}, (\ref{eq:tree_a}).]

\FOR{$\lambda'$ s.t. $\Omega_{\lambda'} \subset \Omega_{\lambda}$}
\IF[$\Omega_{\lambda'}$ is a leaf, {\it i.e.}, (\ref{eq:leaf_a}).]{$\lambda' \in D(h)$ }
\STATE $l=h(\lambda')$.
\STATE $\widetilde{a}_{i,l} = \widetilde{a}_{i,l} + 
|\Omega_{\lambda}|^{-1}|\Omega_{\lambda'}| \alpha_{\gamma,\lambda}$.
\STATE $\widetilde{a}_{i',l} = \widetilde{a}_{i',l} - 
|\Omega_{\lambda}|^{-1}|\Omega_{\lambda'}| \alpha_{\gamma,\lambda}$.
\ELSE
\STATE Continue up to leaves.
\ENDIF
\ENDFOR

\ENDIF
\ENDFOR
\ELSIF[$\Omega_\mu$ is a phantom, {\it i.e.}, (\ref{eq:phantom_a}).]{$\mu \in \bigcup_{r=1}^{N_{\rm R}}  \phan(\Lambda_r)$ }
\FOR{$\lambda \in R^+_F(\gamma)$}
\IF[$\Omega_\lambda$ is a leaf, {\it i.e.}, (\ref{eq:leaf_a}).]{$\lambda \in D(h)$ }
\STATE $l=h(\lambda)$.
\STATE $\widetilde{a}_{i,l} = \widetilde{a}_{i,l} + \alpha_{\gamma,\lambda}$. 

\FOR{$\widehat{\mu}$ s.t. $\Omega_{\mu} \subset \Omega_{\widehat{\mu}}  \wedge
\overline{\Omega_{\gamma}} \cap \overline{\Omega_{\widehat{\mu}}} \neq \varnothing$}
\STATE $i'=h(\widehat{\mu})$.
\STATE $\widetilde{a}_{i',l} = \widetilde{a}_{i',l} - 
|\Omega_{\widehat{\mu}}|^{-1}|\Omega_{\mu}|\alpha_{\gamma,\lambda}$.
\ENDFOR

\ELSIF[$\Omega_\lambda$ is a phantom, {\it i.e.}, (\ref{eq:phantom_a}).]{$\lambda \in \bigcup_{r=1}^{N_{\rm R}}  \phan(\Lambda_r)$ }
\FOR{$\widehat{\lambda}$ s.t. $\widehat{\lambda} \in R_I(\lambda)$}
\IF[$\Omega_{\widehat{\lambda}}$ is a leaf, {\it i.e.}, (\ref{eq:leaf_a}).]{$\widehat{\lambda} \in D(h)$ }
\STATE $l=h(\widehat{\lambda})$.
\STATE $\widetilde{a}_{i,l} = \widetilde{a}_{i,l} + \beta_{\lambda,\widehat{\lambda}}\alpha_{\gamma,\lambda}$.

\FOR{$\widehat{\mu}$ s.t. $\Omega_{\mu} \subset \Omega_{\widehat{\mu}}  \wedge
\overline{\Omega_{\gamma}} \cap \overline{\Omega_{\widehat{\mu}}} \neq \varnothing$}
\STATE $i'=h(\widehat{\mu})$.
\STATE $\widetilde{a}_{i',l} = \widetilde{a}_{i',l} - 
|\Omega_{\widehat{\mu}}|^{-1}|\Omega_{\mu}|\beta_{\lambda,\widehat{\lambda}}\alpha_{\gamma,\lambda}$.
\ENDFOR

\ELSE[$\Omega_{\widehat{\lambda}}$ is within the tree, {\it i.e.}, (\ref{eq:tree_a}).]
\STATE Continue up to leaves.
\ENDIF
\ENDFOR

\ELSE[$\Omega_{\lambda}$ is within the tree, {\it i.e.}, (\ref{eq:tree_a}).]

\STATE Continue up to leaves.

\ENDIF
\ENDFOR

\ELSE[$\Omega_\mu$ is within the tree, {\it i.e.}, (\ref{eq:tree_a}).]

\FOR{$\widehat{\gamma}$ s.t. $\Omega_{\widehat{\gamma}} \subset \Omega_{\gamma}  \wedge
\overline{\Omega_{\widehat{\gamma}}} \cap \overline{\Omega_{\mu}} \neq \varnothing$}

\STATE Current neighbor: 
$\Omega_{\mu'}$ s.t. 
$\Omega_{\mu'} \subset \Omega_{\mu} \wedge
\Gamma^{d'}_{\widehat{\gamma},\mu'}=\overline{\Omega_{\widehat{\gamma}}} \cap \overline{\Omega_{\mu'}}$.
\STATE $\Omega_{\mu'}$ is a leaf: $i'=h(\mu')$.
\FOR{$\lambda \in R^+_F(\widehat{\gamma})$}

\IF[$\Omega_\lambda$ is a leaf, {\it i.e.}, (\ref{eq:leaf_a}).]{$\lambda \in D(h)$ }
\STATE $l=h(\lambda)$.
\STATE $\widetilde{a}_{i,l} = \widetilde{a}_{i,l} 
+ |\Omega_{\gamma}|^{-1}|\Omega_{\widehat{\gamma}}|  \alpha_{\widehat{\gamma},\lambda}$. 
\STATE $\widetilde{a}_{i',l} = \widetilde{a}_{i',l} 
-  \alpha_{\widehat{\gamma},\lambda}$. 

\ELSIF[$\Omega_\lambda$ is a phantom, {\it i.e.}, (\ref{eq:phantom_a}).]{$\lambda \in \bigcup_{r=1}^{N_{\rm R}}  \phan(\Lambda_r)$ }

\FOR{$\widehat{\lambda}$ s.t. $\widehat{\lambda} \in R_I(\lambda)$}

\IF[$\Omega_{\widehat{\lambda}}$ is a leaf, {\it i.e.}, (\ref{eq:leaf_a}).]{$\widehat{\lambda} \in D(h)$ }
\STATE $l=h(\widehat{\lambda})$.
\STATE $\widetilde{a}_{i,l} = \widetilde{a}_{i,l} + 
|\Omega_{\gamma}|^{-1}|\Omega_{\widehat{\gamma}}|
\beta_{\lambda,\widehat{\lambda}}\alpha_{\widehat{\gamma},\lambda}$.
\STATE $\widetilde{a}_{i',l} = \widetilde{a}_{i',l} - 
\beta_{\lambda,\widehat{\lambda}}\alpha_{\widehat{\gamma},\lambda}$.

\ELSE[$\Omega_{\widehat{\lambda}}$ is within the tree, {\it i.e.}, (\ref{eq:tree_a}).]
\STATE Continue up to leaves.

\ENDIF
\ENDFOR

\ELSE[$\Omega_{\lambda}$ is within the tree, {\it i.e.}, (\ref{eq:tree_a}).]

\STATE Continue up to leaves.

\ENDIF

\ENDFOR
\ENDFOR

\ENDIF
\ENDFOR
\ENDFOR
\end{algorithmic}

\section{Photoionization model}
\label{PhotoIonSec}
The photoionization source term $\Sph$ is evaluated using the three-group SP$_3$ model developed in 
\cite{Bourdon:2007} with   Larsen's boundary conditions \cite{Larsen,Liu:2007}.
This model considers $N_g=3$ effective  monochromatic radiative transfer equations.  
As no scattering
of photons is taken into account and since the time scale of photon propagation
is considered short with respect to the streamer propagation,
at each instant of the streamer simulation the photon distribution function $\Psi_{\!l}(\x, \vec\Omega)$
at position $\x$ and direction $\vec\Omega$  fulfills a radiative transfer equation of the form:
\begin{equation}\label{RTE}
{\vec \Omega}\cdot \Dx \Psi_{\!l}(\x, \vec\Omega)+\lambda_l p_{{\rm O}_2} \Psi_{\!l} (\x, \vec\Omega) = 
\frac{1}{4\pi}\frac{p_{\rm q}}{p+p_{\rm q}}\left(\xi\frac{\nu_{\rm u}}{\nu_{\rm i}}\right)\frac{\nu_{\rm i}\ne}{c\,\xi}, 
\quad l=1,\dots,N_g,
\end{equation}
where $l$ indicates discrete wavelengths, $\lambda_l$ is the absorption coefficient, $p_{{\rm O}_2}$ is the partial pressure
of molecular oxygen ($150\,$Torr at atmospheric pressure), $p$ is the total pressure, $p_{\rm q}=30\,$Torr
is the quenching pressure, $\xi=0.1$ is the photoionization efficiency, $\nu_{\rm u}$ is 
the effective excitation coefficient for N$_2$ states responsible for ionizing radiation, 
and $\nu_{\rm i}$ and $\ne$
are, respectively, the previously introduced  ionization coefficient and electron density. 
The term $(\xi\nu_{\rm u}/\nu_{\rm i})$ is
given as a function of the reduced electric field in \cite{Zheleznyak82,Liu:2004}. Finally, $c$ stands for the speed of light.
Let us emphasize that monochromatic equations (\ref{RTE}) have different absorption coefficients
but they all have the same source term 
that depends on the local reduced electric field
$E/N_{\rm air}$, varying therefore in time and space.

The SP$_3$ approximation of (\ref{RTE}) 
leads to a set  of two elliptic equations for functions $\phi_{1,l}(\x)$ and $\phi_{2,l}(\x)$  \cite{Larsen}:
\begin{equation}
\label{SP3}
\left.
\begin{array}{l}
\displaystyle \Dx^2\phi_{1,l}(\x) - \frac{\lambda_l^2\pOtwo^2}{\kappa_1^2}\phi_{1,l}(\x) 
= -\frac{\lambda_l\pOtwo}{\kappa_1^2}
\frac{p_{\rm q}}{p+p_{\rm q}}
\left(\xi\frac{\nu_{\rm u}}{\nu_{\rm i}}\right)\frac{\nu_{\rm i}\ne}{c\,\xi}, \\[1.5ex] 
\displaystyle \Dx^2\phi_{2,l}(\x) - \frac{\lambda_l^2\pOtwo^2}{\kappa_2^2}\phi_{2,l}(\x) = 
-\frac{\lambda_l\pOtwo}{\kappa_2^2}
\frac{p_{\rm q}}{p+p_{\rm q}}
\left(\xi\frac{\nu_{\rm u}}{\nu_{\rm i}}\right)\frac{\nu_{\rm i}\ne}{c\,\xi}, 
\end{array}
\right\}
\end{equation}
with $\kappa_{1,2}=(1/7)(3\pm2\sqrt{6/5})$.
Equations (\ref{SP3}) are coupled through the boundary condition.  
On a boundary surface with neither reflection nor emission, functions $\phi_{1,l}(\x)$ and $\phi_{2,l}(\x)$ 
must verify the following conditions \cite{Larsen,Liu:2007}:
\begin{equation}
\label{SP3BC}
\left.
\begin{array}{l}
\Dx \phi_{1,l}(\x)\cdot \vec n_{s} = 
-\lambda_l\pOtwo\alpha_1  \phi_{1,l}(\x) - \lambda_l\pOtwo\beta_2\phi_{2,l}(\x),\\[1.5ex]
\Dx \phi_{2,l}(\x)\cdot \vec n_{s} = 
-\lambda_l\pOtwo\alpha_2  \phi_{2,l}(\x) - \lambda_l\pOtwo\beta_1\phi_{1,l}(\x),
\end{array}
\right\}
\end{equation}
where $\vec n_{s}$ is the outward unit normal to the boundary surface, $\alpha_{1,2}=(5/96)(34\pm11\sqrt{6/5})$, and 
$\beta_{1,2}=(5/96)(2\pm\sqrt{6/5})$. 
Because $0 < \beta_{1,2}\ll\alpha_{1,2}$ the coupling in  (\ref{SP3}) is weak.
A simple strategy to solve (\ref{SP3}) together with the boundary conditions (\ref{SP3BC}) 
consist in solving the equations independently,
that is, with $\beta_{1,2}=0$
to then iterate and correct the initial approximations
with the inclusion of the $\beta_{1,2}$ coefficients. 
Convergence is attained very rapidly after few iterations (typically three).
The isotropic part of the photon distribution function $\Psi_{\!l}(\x)$
is then written as a linear combination of $\phi_{1,l}(\x)$ and $\phi_{2,l}(\x)$ \cite{Larsen}:
\begin{equation*}\label{PhotonDistr}
\Psi_{\! l}(\x)  = \frac{\gamma_2\phi_{1,l}(\x)-\gamma_1\phi_{2,l}(\x)}{\gamma_2-\gamma_1},
\end{equation*}
with $\gamma_{1,2}=(5/7)(1\pm 3 \sqrt{5/6})$.
The photoionization source term $\Sph(\x)$ can be finally calculated as \cite{Bourdon:2007}:
\begin{equation*}
\Sph(\x) = \sum\limits_{l=1}^{N_{g}} A_l \xi \pOtwo\light \Psi_{\!l}(\x),
\end{equation*}
where  parameters $A_l$ together with $\lambda_l$ are given in Table \ref{TableCoefs}.
\begin{table}[!ht]
\caption{Parameters for three group photoionization model \cite{Bourdon:2007}.}
\label{TableCoefs}
\begin{center}
\begin{tabular}{ccc} 
\toprule
$l$ & $A_l\, [{\rm cm}^{-1}\, {\rm Torr}^{-1}]$ & $\lambda_l [{\rm cm}^{-1}\, {\rm Torr}^{-1}]$ \cr
\midrule
1 & 0.0067 &  0.0447 \cr
2 & 0.0346 &  0.1121 \cr
3 & 0.3059 &  0.5994 \cr
\bottomrule
\end{tabular}
\end{center}
\end{table}

\bibliographystyle{plain}
 \bibliography{biblio_Poisson.bib}

\end{document}